\documentclass[a4paper,11pt]{article}

\usepackage{amsmath}
\usepackage{booktabs}
\usepackage{graphicx}
\usepackage[english]{babel}
\usepackage{nicefrac}
\usepackage{subcaption}

\title{An integrated mixed integer program model for the two-level capacitated vehicle routing problem with extensions}

\author{Congzheng Liu\thanks{Decision Lab Ltd,
London SE11 5JH, UK.
Email: {\tt \{joshua.liu@decisionlab.co.uk\}}}
} 

\begin{document}
\maketitle

\begin{abstract}
This paper introduces the two-level capacitated vehicle routing problem (2S-CVRP). This problem combines the two-level bin packing problem and the vehicle routing problem into an integrated framework. The problem itself is an NP-hard problem and it can be seen as an extension to the traditional capacitated vehicle routing problem (CVRP). We propose this extension as it enable one to model more real-life applications in logistics. A mixed integer program (MIP) model is presented for the problem. Our MIP model includes an extensive set of constraints encountered in real-world applications. The validity of the model is tested on both artificial and real-life instances. 
\\*[2mm]
\end{abstract}


\section{Introduction}

This paper introduces and studies an \emph{Operation Research} (OR) problem called the two-level capacitated vehicle routing problem (2S-CVRP). It is a NP-hard problem with real-world applications in logistics. This problem states that boxes with certain measurements and pre-determined destinations are to be packed inside available pallets without exceeding their capacities. These pallets then need to be assigned to trucks, and the trucks are routed in order to deliver all boxes to their destination. The goal is to minimise the overall costs. The major difficulty of 2S-CVRP is that each pallet is not necessarily delivered to only one destination.

In early works on the field of OR, the bin packing problem (BPP) and the vehicle routing problem (VRP) are usually treated independently. Great amount of works have been done in the respective fields. Specifically, the BPP considers a situation where boxes of different sizes must be packed into a finite number of pallets, each of a fixed given capacity, in a way that minimises the cost of pallets used. The two-level BPP, as an extension, introduces an additional level, where the pallets must be fitted into trucks. It can be viewed as a variant that considers additional setting to improve efficiency in packing. The earliest paper that developed the bin packing model, to our knowledge, is \cite{Ka60} (Translation of a 1939 paper in Russian). Following that, the concept of column generation, an efficient algorithm for solving large linear programs, was proposed in \cite{GG61}. Since then, the literature on the bin packing problem is vast (see e.g., \cite{SP92,DIM16} for review).

Vehicle routing problem (VRP), on the other hand, is to find the optimal routes for multiple vehicles visiting a set of locations, where the ``optimal routes" can be defined as the routes with the least total cost, for instance. The capacitated VRP, as an extension, takes the capacity constraints of vehicles into account (for example, see \cite{ISV07,GILM08,WZZL15}). The idea of VRP first appeared in a paper by Dantzig and Ramser in 1959, in which an algorithmic approach was written and was applied to petrol deliveries. Later on, in \cite{CW64}, the Dantzig and Ramser's approach was improved using an effective greedy algorithm. Common exact approaches to model the VRP include: branch and bound approach, cutting plane approach, network flow approach and dynamic program (\cite{Zi08}).

The idea of the two-level capacitated vehicle routing problem (2S-CVRP) is to combine the two-level BPP and the capacitated VRP into an integrated framework. This is motivated by the fact that in real world distribution problem, the ``packing" and ``transport" decision are usually required to be considered simultaneously. Typical examples include groceries delivery, parcel postage and air cargo loading. Such distribution problems often have similar characteristics: i). The packing stage involves more than one level, e.g., pallet level, truck level. ii). The outcome of the packing stage will influence the decision on the routing stage. iii). The overall cost is more of interest than the independent costs.

Our aim is to provide a rich and realistic mathematical representation of this 2S-CVRP. The structure of this paper will be organised as follows. In Section \ref{se:problem}, we briefly describe the objective and constraints of the 2S-CVRP. In Section \ref{se:model}, we present a mixed integer program (MIP) model to the 2S-CVRP. Section \ref{se:extension} gives extensions to the model. The model is tested on both artificial and real-life instances in Section \ref{se:experiment}. Finally, Section \ref{se:conclusion} concludes the paper.


\section{Problem definition}
\label{se:problem}

Given a set $I$ of boxes, a set $J$ of pallets, a set $K$ of trucks and a set $D$ of destinations, the 2S-CVRP consists of delivering all boxes to their pre-determined destination using pallets and trucks with minimum costs. In particular, all boxes have to be packed into pallets first in order to be loaded onto trucks. Direct packing from boxes to trucks is forbidden under all circumstances. Every box $i \in I$ is characterised by a volume $v_i$, every pallet $j \in J$ is characterised by a capacity $V_j$ and a fix cost $c_j$. In general, we do not require the value of $V_j$ and $c_j$ to be identical for all $j$. However, for computational simplicity, we assume the outer volume of each pallet is equal to its inner capacity. Additionally, we let $T_k$ and $C_k$ denote the capacity and fix cost of truck $k \in K$, respectively. For $i \in I, d\in D$, we let $g^0_{id}$ takes value 1 if box $i$ needs to be delivered to destination $d$, and takes value 0 otherwise. Finally, for $d,d' \in \{D \cap d_0\}$, we let $c'_{dd'}$ denote the variable cost of travelling from destination $d$ to destination $d'$, where $d_0$ denotes the depot.

To ensure at least one feasible solution, some conditions are assumed to be satisfied: for example, the volume of each box is supposed to be less or equal to the maximum capacity of the pallets ($v_i \leq \max_j V_j, \text{for} \; i \in I$). This is also true for the relationship between pallets and trucks ($V_i \leq \max_k T_k, \text{for} \; j \in J$). We also make sure that each box is to be delivered to only one destination. We should note, however, each pallet is not necessarily delivered to only one destination, which is the major difference between 2S-CVRP and classic capacitated VRP. The fundamental of 2S-CVRP is to consider the following problem:
\[
\begin{aligned}
\min \quad & \text{total costs of packing and routing}\\
\textrm{s.t.} \quad & \text{each box assigned to exactly one pallet}\\
& \text{each pallet assigned to exactly one truck}\\
& \text{the total volume of the boxes inside a pallet} \leq \text{pallet capacity}\\
& \text{the total volume of the pallets inside a truck} \leq \text{truck capacity}\\
& \text{each box delivered to its destination}
\end{aligned}
\]


\section{Mixed integer program model}
\label{se:model}

In this section, we consider the mathematically formulation of the 2S-CVRP under a mixed integer program (MIP). In Subsection \ref{sub:variable}, we introduce the most basic variables. Subsection \ref{sub:objective} considers the objective function of 2S-CVRP. Then, we explain all constraints one by one in Subsection \ref{sub:constraints}.

\subsection{Variables}
\label{sub:variable}

Let's begin with the most basic variables that concerns the packing process itself. For now, we assume the packing is one dimensional. Multi-dimensional packing will be introduced as extensions in Section \ref{se:extension}. Here are the basic variables to used in the model. Note that the subscripts relate to indices and superscripts relate to levels.
\[
\begin{aligned}
& p^0_{ij} = \begin{cases}
1, & \text{if box} \; i \; \text{is in pallet} \; j, \\
0, & \text{otherwise,}
\end{cases} & \text{for} \; i \in I, j \in J\\
& u^0_j = \begin{cases}
1, & \text{if pallet} \; j \; \text{is used,}\\
0, & \text{otherwise,}
\end{cases} & \text{for} \; j \in J\\
& g^1_{jd} = \begin{cases}
1, & \text{if pallet} \; j \; \text{is calling at} \; d,\\
0, & \text{otherwise,}
\end{cases} & \text{for} \; j \in J, d \in D\\
& p^1_{jk} = \begin{cases}
1, & \text{if pallet} \; j \; \text{is in truck} \; k, \\
0, & \text{otherwise,}
\end{cases} & \text{for} \; j \in J, k \in K\\
& u^1_k = \begin{cases}
1, & \text{if truck} \; k \; \text{is used,}\\
0, & \text{otherwise,}
\end{cases} & \text{for} \; k \in K\\
& \eta_{jkd}^1 = \begin{cases}
1, & \text{if pallet} \; j \; \text{in} \; k \; \text{is delivered to} \; d,\\
0, & \text{otherwise,}
\end{cases} & \text{for} \; j \in J,k \in K,d \in D\\
& g^2_{kd} = \begin{cases}
1, & \text{if truck} \; k \; \text{is calling at} \; d,\\
0, & \text{otherwise,}
\end{cases} & \text{for} \; k \in K, d \in D\\
& g^s_{kdd'} = \begin{cases}
1, & \text{if truck} \; k \; \text{uses the route} \; d \; \text{to} \; d',\\
0, & \text{otherwise,}
\end{cases} & \text{for} \; k \in K, d,d' \in \{D \cap d_0\}\\
\end{aligned}
\]
Here we use $d_0$ to denote the depot. Each truck in use needs to leave from the depot to make delivery and come back to depot once the delivery is finished.

\subsection{Objective function}
\label{sub:objective}

The objective function consists in minimising the summation of the fix cost of pallet, the fix cost of truck and the variable cost of truck calling at destinations:
\begin{equation}
    \min_{u^0_j, u^1_k, g^s_{kdd'}} \quad \sum_{j \in J} u^0_j c_j + \sum_{k \in K} u^1_k C_k + \sum_{k \in K, d,d' \in \{D \cap d_0\}} g^s_{kdd'} c'_{dd'}
\end{equation}
The first two parts of the objective function denote the packing cost of pallet and truck, identical to the ones in bin packing problem, while the third part denotes the routing cost, which is very similar to the objective of classical vehicle routing problem.

\subsection{Constraints}
\label{sub:constraints}

First, we need to make sure that each box/pallet assigned to exactly one pallet/truck and does not exceed capacities.
\begin{eqnarray}
\label{eq:cap_0}
& \displaystyle \sum_{i \in I} v_i p^0_{ij} \leq V_j u^0_j & \text{for} \; j \in J\\
\label{eq:assign_0}
& \displaystyle \sum_{j \in J} p^0_{ij} = 1 & \text{for} \; i \in I\\
\label{eq:cap_1}
& \displaystyle \sum_{j \in J} V_j p^1_{jk} \leq T_k u^1_k & \text{for} \; k \in K\\
\label{eq:assign_1}
& \displaystyle \sum_{k \in K} p^1_{jk} = u_j^0 & \text{for} \; j \in J
\end{eqnarray}
The maximum capacity of each pallet $j$ cannot be exceeded, which is ensured by constraints \eqref{eq:cap_0}. Constraints \eqref{eq:assign_0} verify that each box is allocated to exactly one pallet. The relationship between pallet and truck is defined similarly by constraints \eqref{eq:cap_1} - \eqref{eq:assign_1}. Then, we consider the destination assignment:
\begin{eqnarray}
\label{eq:des_1}
& \displaystyle \sum_{i \in I} p_{ij}^0 g_{id}^0 \leq g^1_{jd} \, \Upsilon & \text{for} \; j \in J, d \in D\\
\label{eq:des_2}
& 2 \, \eta_{jkd}^1 \leq p_{jk}^1 + g_{jd}^1 \leq \eta_{jkd}^1 +1 & \text{for} \; j \in J,k \in K, d \in D\\
\label{eq:des_3}
& \displaystyle \sum_{j \in J} \eta_{jkd}^1 \leq g^2_{kd} \, \Upsilon & \text{for} \; k \in K, d \in D
\end{eqnarray}
Constraints \eqref{eq:des_1} make sure that each pallet is calling at a destination if it contains boxes need to be delivered to that destination. Then, constraints \eqref{eq:des_2} - \eqref{eq:des_3} ensure that each truck is calling at a destination if it contains pallets calling at that destination. This is done by adding an auxiliary variable $\eta_{jkd}^1$ that takes value 1 if $p_{jk}^1$ and $g_{jd}^1$ are both 1, and takes value 0 otherwise.

The following constraints ensure that the trucks are assigned to destinations accordingly:
\begin{eqnarray}
\label{eq:routeallow_1}
& 2 \, g^s_{kdd'} \leq g^2_{kd} + g^2_{kd'} & \text{for} \; k \in K,d,d' \in D\\
\label{eq:routeallow_2}
& g^s_{kd_0d} + g^s_{kdd_0} \leq 2 \,g^2_{kd} & \text{for} \; k \in K,d \in D\\
\label{eq:routeassign_1}
& \sum_{d \in \{D \cap d_0\}} g^s_{kdd'} = g^2_{kd'} & \text{for} \; k \in K,d' \in \{D \cap d_0\}\\
\label{eq:routeassign_2}
& \sum_{d' \in \{D \cap d_0\}} g^s_{kdd'} = g^2_{kd} & \text{for} \; k \in K,d \in \{D \cap d_0\}
\end{eqnarray}
Constraints \eqref{eq:routeallow_1} - \eqref{eq:routeallow_2} describe that the route $d$ to $d'$ (or the depot $d_0$) can only be used by truck $k$ if the truck is calling at both destination $d$ and destination $d'$ (or the depot $d_0$). Constraints \eqref{eq:routeassign_1} - \eqref{eq:routeassign_2} make sure that each destination is entered once and left once by each truck if and only if the truck is calling at that destination. 

To eliminate the possibility of subtour. We introduce an auxiliary variable following \cite{MTZ60}:
\[
\begin{aligned}
e_{kd} & \quad \text{integer auxiliary variable at} \; d \; \text{for truck} \; k & \text{for} \; k \in K, d \in D\\
\end{aligned}
\]
Subtour elimination constraints are given as:
\begin{eqnarray}
\label{eq:detour_1}
& e_{kd} - e_{kd'} + |D| g^s_{kdd'} \leq |D| -1 & \text{for} \; k \in K,d,d' \in D\\
\label{eq:detour_2}
&  1 \leq e_{kd} \leq |D| & \text{for} \; k \in K,d \in D
\end{eqnarray}
Constraints \eqref{eq:detour_1} - \eqref{eq:detour_2} make sure all subtours are eliminated. If truck $k$ drives from $d$ to $d'$, $g^s_{kdd'}=1$ and constraints \eqref{eq:detour_1} become $e_{kd'} \geq e_{kd} +1$. This makes sure that $e_{kd'}$ is always greater than $e_{kd}$, forcing the auxiliary variable $e$ to increase at the arrival of each new destination. In this case, the truck cannot revisit any previous destination. If truck $k$ does not drive from $d$ to $d'$, the constraints \eqref{eq:detour_1} become $e_{kd'} - 1 \geq e_{kd} -|D|$, which is always valid with the help of constraints \eqref{eq:detour_2}. 


\section{Multi-dimensional packing}
\label{se:extension}

In this section, we extend the 2S-CVRP model with a more complex case, where the packing capacity of the pallet/truck is multi-dimensional. For the packing at the pallet level, we consider it to be 3-dimensional. This will be discussed in Subsection \ref{sub:pallet3d}. For the packing at the truck level, we consider it to be 2-dimensional, since most trucks are one storey. This will be discussed in Subsection \ref{sub:truck2d}. 

\subsection{Pallet packing}
\label{sub:pallet3d}
The 3-dimensional pallet packing considers more geometry constraints, including orthogonal placement, no overlap, orientation constraints. Moreover, we would like to make sure the boxes are packed considering the delivery sequence. For instance, the box to be delivered at a ``later'' destination on the route should not be placed on the top of the box to be delivered at an ``early'' destination, if they are packed in the same pallet. For these purposes, we introduce some new parameters and variables:

\[
\begin{aligned}
& l^0_i \times w^0_i \times h^0_i  \quad \text{Length} \times \text{width} \times \text{height of box} \; i, \qquad \quad & \text{for} \; i \in I\\
& l^1_j \times w^1_j \times h^1_j  \quad \text{Length} \times \text{width} \times \text{height of pallet} \; j, \qquad \quad & \text{for} \; j \in J\\
& (x^0_i,y^0_i,z^0_i) \qquad \text{location of the front left bottom of box} \; i, & \text{for} \; i \in I\\
& (\hat{x}^0_i,\hat{y}^0_i,\hat{z}^0_i) \qquad \text{location of the rear right top of box} \; i, & \text{for} \; i \in I\\
\end{aligned}
\]
\[
\begin{aligned}
& r^0_{iab} = \begin{cases}
1, & \text{if side} \; b \; \text{of box} \; i \; \text{is along the} \; a \text{-axis of pallet,}\\
0, & \text{otherwise,}
\end{cases} & \text{for} \; i \in I\\
& x_{ii'}^{p0} = \begin{cases}
1, & \text{if box} \; i \; \text{is on the right of box} \; i' \; (\hat{x}^0_{i'} \leq x^0_i),\\
0, & \text{otherwise} \; (x^0_i < \hat{x}^0_{i'}),
\end{cases} & \text{for} \; i,i' \in I\\
& y_{ii'}^{p0} = \begin{cases}
1, & \text{if box} \; i \; \text{is behind box} \; i' \; (\hat{y}^0_{i'} \leq y^0_i),\\
0, & \text{otherwise} \; (y^0_i < \hat{y}^0_{i'}),
\end{cases} & \text{for} \; i,i' \in I\\
& z_{ii'}^{p0} = \begin{cases}
1, & \text{if box} \; i \; \text{is above box} \; i' \; (\hat{z}^0_{i'} \leq z^0_i),\\
0, & \text{otherwise} \; (z^0_i < \hat{z}^0_{i'}),
\end{cases} & \text{for} \; i,i' \in I\\
& p_{ik}^g = 
\begin{cases}
1, & \text{if box} \; i \; \text{is in truck} \; k,\\
0, & \text{otherwise,}
\end{cases} & \text{for} \; i \in I,k \in K\\
& \theta_{ii'}^0 = \begin{cases}
1, & \text{if box} \; i' \; \text{is delivered immediately after} \; i,\\
0, & \text{otherwise,}
\end{cases} & \text{for} \; i,i' \in I\\
\end{aligned}
\]
We have $a,b \in \{1,2,3\}$. The index $b$ indicates the side of the box, that is $b \in \{l^0 := 1, w^0 := 2, h^0 := 3\}$, whereas $a$ indicates the axis, that is $a \in \{x := 1, y := 2, z := 3\}$. Without loss of generality, the axes of
the coordinate system are assumed to be placed so that the length $l^1_j$ (resp. width $w^1_j$, height
$h^1_j$) of the pallet j lies on the x-axis (resp. y-axis, z-axis) $\text{for} \; j \in J$. The origin of this
coordinate system lies on the front left bottom corner of the pallets.

First, we consider some basic geometric constraints concerning orthogonal placement:
\begin{eqnarray}
\label{eq:size_0_x}
& \hat{x}^0_i \leq \displaystyle \sum_{j \in J} l^1_j p^0_{ij} & \text{for} \; i \in I\\
\label{eq:size_0_y}
& \hat{y}^0_i \leq \displaystyle \sum_{j \in J} w^1_j p^0_{ij} & \text{for} \; i \in I\\
\label{eq:size_0_z}
& \hat{z}^0_i \leq \displaystyle \sum_{j \in J} h^1_j p^0_{ij} & \text{for} \; i \in I
\end{eqnarray}
Constraints \eqref{eq:size_0_x} - \eqref{eq:size_0_z} ensure that the boxes do not exceed their pallet size.

For the orientation, constraints \eqref{eq:rot_0_x} - \eqref{eq:rot_0_b} describe that the boxes can rotate orthogonally in the pallet. 

\begin{eqnarray}
\label{eq:rot_0_x}
& \hat{x}^0_i - x^0_i = r^0_{i11}l^0_i + r^0_{i12}w^0_i + r^0_{i13}h^0_i & \text{for} \; i \in I\\
\label{eq:rot_0_y}
& \hat{y}^0_i - y^0_i = r^0_{i21}l^0_i + r^0_{i22}w^0_i + r^0_{i23}h^0_i & \text{for} \; i \in I\\
\label{eq:rot_0_z}
& \hat{z}^0_i - z^0_i = r^0_{i31}l^0_i + r^0_{i32}w^0_i + r^0_{i33}h^0_i & \text{for} \; i \in I
\end{eqnarray}
\begin{eqnarray}
\label{eq:rot_0_a}
& \displaystyle \sum_{a \in \{1,2,3\}} r^0_{iab} = 1 & \text{for} \; i \in I, b \in \{1,2,3\}\\
\label{eq:rot_0_b}
& \displaystyle \sum_{b \in \{1,2,3\}} r^0_{iab} = 1 & \text{for} \; i \in I, a \in \{1,2,3\}
\end{eqnarray}

The following constraints ensure that there is no overlap:
\begin{eqnarray}
\label{eq:pos_0}
& x_{ii'}^{p0} + x_{i'i}^{p0} + y_{ii'}^{p0} +  y_{i'i}^{p0} + z_{ii'}^{p0} + z_{i'i}^{p0} \geq p^0_{ij} + p^0_{i'j} - 1 & \text{for} \; i,i' \in I, j \in J \quad\\
\label{eq:pos_0_x_1}
& \hat{x}^0_{i'} \leq x^0_i + (1 - x^{p0}_{ii'}) \Upsilon & \text{for} \; i,i' \in I\\
\label{eq:pos_0_x_2}
& x^0_i + \epsilon \leq \hat{x}^0_{i'} + x^{p0}_{ii'} \Upsilon & \text{for} \; i,i' \in I\\
\label{eq:pos_0_y_1}
& \hat{y}^0_{i'} \leq y^0_i + (1 - y^{p0}_{ii'}) \Upsilon & \text{for} \; i,i' \in I\\
\label{eq:pos_0_y_2}
& y^0_i + \epsilon \leq \hat{y}^0_{i'} + y^{p0}_{ii'} \Upsilon & \text{for} \; i,i' \in I\\
\label{eq:pos_0_z_1}
& \hat{z}^0_{i'} \leq z^0_i + (1 - z^{p0}_{ii'}) \Upsilon & \text{for} \; i,i' \in I
\end{eqnarray}
If the variables $x_{ii'}^{p0}$, $x_{i'i}^{p0}$, $y_{ii'}^{p0}$, $y_{i'i}^{p0}$, $z_{ii'}^{p0}$ or $z_{i'i}^{p0}$ equal 1, the two boxes $i$ and $i'$ do not overlap along any of the axes. To prevent having two boxes occupying a same portion of space, it is sufficient
to allow no overlap along at least one of the axes, that is, at least one of these variables must equal 1. It leads to constraints \eqref{eq:pos_0}. Constraints \eqref{eq:pos_0_x_1} - \eqref{eq:pos_0_z_1} ensure that $x_{ii'}^{p0} = 1$ (resp.
$y_{ii'}^{p0} =1$, $z_{ii'}^{p0} =1$), if and only if $x^0_i \geq \hat{x}^0_{i'}$ (resp.
$y^0_i \geq \hat{y}^0_{i'}$, $z^0_i \geq \hat{z}^0_{i'}$).

Finally, we consider our preference of the box allocation.
\begin{eqnarray}
\label{eq:allocation_1}
& 2 \, p_{ik}^g \leq p^0_{ij} + p^1_{jk} \leq p_{ik}^g +1 & \text{for} \; i \in I,j \in J,k \in K\\
\label{eq:allocation_2}
& 4 \, \theta_{ii'}^0 \leq g^s_{kdd'} + g_{id}^0 g_{i'd'}^0 + p^g_{ik} + p^g_{i'k} & \text{for} \; i,i' \in I,k \in K,d,d' \in D \qquad\\
\label{eq:allocation_3}
& g^s_{kdd'} + g_{id}^0 g_{i'd'}^0 + p^g_{ik} + p^g_{i'k} \leq \theta_{ii'}^0 +3 & \text{for} \; i,i' \in J,k \in K,d,d' \in D \qquad\\
\label{eq:allocation_4}
& \theta_{ii'}^0 + p^0_{ij} + p^0_{i'j} \leq 3\, (1-z^{p0}_{i'i}) & \text{for} \; i,i' \in I, j \in J
\end{eqnarray}
Constraints \eqref{eq:allocation_1} define new variable $p_{ik}^g$ that takes value 1 if box $i$ is in truck $k$ (both $p_{ij}^0$ and $p_{jk}^1$ take value 1), and takes value 0 otherwise. Constraints \eqref{eq:allocation_2} - \eqref{eq:allocation_3} define auxiliary variable $\theta_{ii'}^0$, that takes value 1 if and only if all four conditions are met. Lastly, constraint \eqref{eq:allocation_4} ensure that box $i'$ cannot be put on the top of box $i$ if box $i'$ is delivered immediately after box $i$ and in the same pallet (and truck).

\subsection{Truck packing}
\label{sub:truck2d}

For the packing at the truck level, we consider that: i). The geometric constraints must be met. ii). The sequence of loading the truck is in reverse order of the delivery sequence. For simplicity, we consider the problem to be 2-dimensional, since most trucks are one storey. Moreover, the pallets are allowed to be rotated horizontally, but not vertically. For this propose, we introduce some new parameters and variables:
\[
\begin{aligned}
& l^2_k \times w^2_k \quad \text{Length} \times \text{width of truck} \; k, \qquad \quad & \text{for} \; k \in K\\
& (x^1_j,y^1_j) \qquad \text{location of the front left corner of} \; j, & \text{for} \; j \in J\\
& (\hat{x}^1_j,\hat{y}^1_j) \qquad \text{location of the rear right corner of} \; j, & \text{for} \; j \in J\\
& r^1_{jab} = \begin{cases}
1, & \text{if side} \; b \; \text{of pallet} \; j \; \text{is along the} \; a \text{-axis,}\\
0, & \text{otherwise,}
\end{cases} & \text{for} \; j \in J\\
& x_{jj'}^{p1} = \begin{cases}
1, & \text{if pallet} \; j \; \text{is on the right of pallet} \; j',\\
0, & \text{otherwise},
\end{cases} & \text{for} \; j,j' \in J\\
& y_{jj'}^{p1} = \begin{cases}
1, & \text{if pallet} \; j \; \text{is behind pallet} \; j' \; ,\\
0, & \text{otherwise},
\end{cases} & \text{for} \; j,j' \in J\\
& g_{jj'}^{dd'} = \begin{cases}
1, & \text{if pallet} \; j \; \text{calling} \; d \; \text{and} \; j' \; \text{calling} \; d',\\
0, & \text{otherwise},
\end{cases} & \text{for} \; j,j' \in J,d,d' \in D\\
& \theta_{jj'}^1 = \begin{cases}
1, & \text{if} \; j' \; \text{is delivered immediately after} \; j,\\
0, & \text{otherwise,}
\end{cases} & \text{for} \; j,j' \in J\\
\end{aligned}
\]

Without loss of generality, we assume the axes of the coordinate system to be placed so that the length $l^2_j$ and width $w^2_j$ of the truck k lies on the x-axis and y-axis $\text{for} \; k \in K$. The origin of this
coordinate system lies on the front left corner of the trucks. Therefore, we would like to ensure that the pallet to be delivered to a ``later'' destination cannot be allocated to the left of the pallet to be delivered to an ``early'' destination.

Here are the constraints for the truck packing:

\begin{eqnarray}
\label{eq:size_1_x}
& \hat{x}^1_j \leq \displaystyle (1 - u_j^0) \Upsilon + \sum_{k \in K} l^2_k p^1_{jk} & \text{for} \; j \in J\\
\label{eq:size_1_y}
& \hat{y}^1_j \leq \displaystyle (1 - u_j^0) \Upsilon + \sum_{k \in K} w^2_k p^1_{jk} & \text{for} \; j \in J\\
\label{eq:rot_1_x}
& \hat{x}^1_j - x^1_j = r^1_{j11}l^1_j + r^1_{j12}w^1_j & \text{for} \; j \in J\\
\label{eq:rot_1_y}
& \hat{y}^1_j - y^1_j = r^1_{j21}l^1_j + r^1_{j22}w^1_j & \text{for} \; j \in J\\
\label{eq:rot_1_a}
& \displaystyle \sum_{a \in \{1,2\}} r^1_{jab} = 1 & \text{for} \; j \in J, b \in \{1,2\}\\
\label{eq:rot_1_b}
& \displaystyle \sum_{b \in \{1,2\}} r^1_{jab} = 1 & \text{for} \; j \in J, a \in \{1,2\}\\
\label{eq:pos_1}
& x_{jj'}^{p1} + x_{j'j}^{p1} + y_{jj'}^{p1} +  y_{j'j}^{p1} \geq p^1_{jk} + p^1_{j'k} - 1 & \text{for} \; j,j' \in J, k \in K
\end{eqnarray}
\begin{eqnarray}
\label{eq:pos_1_x_1}
& \hat{x}^1_{j'} \leq x^1_j + (1 - x^{p1}_{jj'}) \Upsilon & \text{for} \; j,j' \in J\\
\label{eq:pos_1_x_2}
& x^1_j + \epsilon \leq \hat{x}^1_{j'} + x^{p1}_{jj'} \Upsilon & \text{for} \; j,j' \in J\\
\label{eq:pos_1_y_1}
& \hat{y}^1_{j'} \leq y^1_j + (1 - y^{p1}_{jj'}) \Upsilon & \text{for} \; j,j' \in J\\
\label{eq:calling_1}
& 2 \, g_{jj'}^{dd'} \leq g^1_{jd} + g^1_{j'd'} \leq g_{jj'}^{dd'} +1 & \text{for} \; j,j' \in J,d,d' \in D\\
\label{eq:calling_2}
& 4 \, \theta_{jj'}^1 \leq g^s_{kdd'} + g_{jj'}^{dd'} + p^1_{jk} + p^1_{j'k} & \text{for} \; j,j' \in J,k \in K,d,d' \in D \qquad\\
\label{eq:calling_3}
& g^s_{kdd'} + g_{jj'}^{dd'} + p^1_{jk} + p^1_{j'k} \leq \theta_{jj'}^1 +3 & \text{for} \; j,j' \in J,k \in K,d,d' \in D \quad\\
\label{eq:order}
& \theta_{jj'}^1 \leq 1-x^{p1}_{jj'} & \text{for} \; j,j' \in J
\end{eqnarray}
Constraints \eqref{eq:size_1_x} - \eqref{eq:pos_1_y_1} are similar to constraints \eqref{eq:size_0_x} - \eqref{eq:pos_0_z_1} in the pallet packing level, in addition to a minor adjustment in order to prevent infeasibility. Constraints \eqref{eq:calling_1} define new variable $g_{jj'}^{dd'}$ that takes value 1 if pallet $j$ is calling at destination $d$ and pallet $j'$ is calling at destination $d'$, and takes value 0 otherwise. Constraints \eqref{eq:calling_2} - \eqref{eq:calling_3} work similarly, but with four conditions. Lastly, constraint \eqref{eq:order} ensure that pallet $j$ cannot be put on the right of pallet $j'$ if pallet $j'$ is delivered immediately after pallet $j$ and on the same truck.


\section{Computational experiments}
\label{se:experiment}

In this section, we perform computational experiments to check the validity of the model. In Subsection \ref{sub:1D}, we consider the 2S-CVRP problem with one-dimensional packing, and \ref{sub:3D}, we consider the case with multi-dimensional packing. Finally, in Subsection \ref{sub:real}, we examine the 2S-CVRP on a real-life instance. All the experiments were performed on an Apple M1 Pro (2021). All the MIPs were solved using the {\tt Gurobi} solver embedded in {\tt Pyomo} in Python. 

\subsection{One-dimensional}
\label{sub:1D}

In order to test the proposed 2S-CVRP model, we adopt 7 groups of randomly generated instances:
\begin{itemize}
    \item ins-1: Set: $|I| = 5$, $|J| = 5$, $|K| = 1$, $|D| = 5$; Parameter: $v \in [1,5]$, $V \in [5,15]$, $T \in [25,28]$, $c \in [5,10]$, $C \in [25,30]$ , $c' \in [0,15]$.
    \item ins-2: Set: $|I| = 5$, $|J| = 5$, $|K| = 2$, $|D| = 5$; Parameter: $v \in [1,8]$, $V \in [8,15]$, $T \in [20,25]$, $c \in [5,10]$, $C = [25,35]$ , $c' \in [0,15]$.
    \item ins-3: Set: $|I| = 9$, $|J| = 5$, $|K| = 3$, $|D| = 5$; Parameter: $v \in [1,6]$, $V \in [10,15]$, $T \in [20,25]$, $c \in [5,10]$, $C = [25,35]$ , $c' \in [0,15]$.
    \item ins-4: Set: $|I| = 9$, $|J| = 5$, $|K| = 3$, $|D| = 5$; Parameter: $v \in [1,5]$, $V \in [8,12]$, $T \in [15,20]$, $c \in [3,10]$, $C = [25,35]$ , $c' \in [0,15]$.
    \item ins-5: Set: $|I| = 11$, $|J| = 6$, $|K| = 3$, $|D| = 5$; Parameter: $v \in [1,5]$, $V \in [8,15]$, $T \in [15,25]$, $c \in [3,12]$, $C = [25,40]$ , $c' \in [0,15]$.
    \item ins-6: Set: $|I| = 20$, $|J| = 10$, $|K| = 6$, $|D| = 5$; Parameter: $v \in [1,4]$, $V \in [5,15]$, $T \in [15,25]$, $c \in [2,5]$, $C = [3,6]$ , $c' \in [0,15]$.
    \item ins-large: Set: $|I| = 20$, $|J| = 10$, $|K| = 6$, $|D| = 16$; Parameter: $v \in [1,4]$, $V \in [5,15]$, $T \in [15,25]$, $c \in [2,5]$, $C = [3,6]$ , $c' \in [0,13]$.
\end{itemize}
Using as benchmark, two traditional disjoint models are considered, both of which solve the two-level bin packing problem and the vehicle routing problem individually. 
\begin{itemize}
    \item 1DBP+VRP: The model generates the packing procedure with optimal cost in the first phase. Then, designs the optimal delivery routes in the second phase.
    \item 1DVRP+BP: The model groups boxes according their destinations in the first phase. Then, packs and routes in the second phase.
\end{itemize}
Since the problem is assumed to be one-dimensional, we denote the proposed integrated method as 2S-1D. The results can be seen in Table \ref{tab:2S-1Dvs1DBP+1DVRP}.

\begin{table}[htb]
\caption{2S-1D vs. 1DBP+VRP/1DVRP+BP cost \& time (in seconds)}
\label{tab:2S-1Dvs1DBP+1DVRP}
\centering
\resizebox{\linewidth}{!}{
\begin{tabular}{ccccccc}
\toprule
\multicolumn{1}{c}{} & \multicolumn{2}{c}{\textbf{Integrated method}} & \multicolumn{4}{c}{\textbf{Disjoint method}}\\
\cmidrule(l{3pt}r{3pt}){2-3} \cmidrule(l{3pt}r{3pt}){4-7}
& 2S-1D cost & 2S-1D time & 1DBP+VRP cost & 1DBP+VRP time & 1DVRP+BP cost & 1DVRP+BP time\\
\midrule
ins-1 & 61 & 0.04 & 61 & 0.19 & 61 & 0.20\\
ins-2 & 114 & 0.08 & 121 & 0.12 & 118 & 0.19\\
ins-3 & 122 & 0.31 & 142 & 1.30 & 179 & 0.30\\
ins-4 & 109 & 0.27 & 122 & 1.38 & 124 & 1.40\\
ins-5 & 114 & 0.51 & 133 & 1.37 & 129 & 1.32\\
ins-6 & 57 & 25.62 & 101 & 20.91 & infeasible & infeasible\\
ins-large & 58 & 69.71 & 73 & 45.27 & 68 & 44.23\\
\bottomrule
\end{tabular}}
\end{table}

It can be seen from Table \ref{tab:2S-1Dvs1DBP+1DVRP} that the integrated 2S-1D method can outperform both the disjoint 1DBP+VRP method and 1DVRP+BP method in regard to the cost in most cases, while requiring similar running time. This is due to the fact that the 2S-1D method optimises the overall cost directly, but both the 1DBP+VRP method and the 1DVRP+BP method treat the bin packing cost and the routing cost individually. In other word, the optimal solution of the disjoint methods will always be upper bounds of the 2S-1D optimal solution. We breakdown the case with instance `ins-6' to observe the detailed differences (values in brackets denote [capacity; cost] for pallets/trucks and denote [cost] for routes):
\begin{itemize}
    \item 2S-1D: Pallet used: J1[15,5], J2[9,4], J3[8,3], J4[8,3], J5[8,3], J10[3,1]; Truck used: K1[25,6], K2[20,5], K3[12,3]; Route used: (D0, D3, D1, D0)[7], (D0, D2, D4, D5, D0)[13], (D0, D2, D0)[4].
    \item 1DBP+VRP: Pallet used: J1[15,5], J2[9,4], J3[8,3], J5[8,3], J6[4,2], J10[3,1]; Truck used: K1[25,6], K3[12,3], K5[12,3]; Route used: (D0, D3, D1, D2, D4, D5, D0)[19], (D0, D3, D1, D4, D5, D0)[26]. (D0, D5, D4, D1, D3, D0)[26].
\end{itemize}

In the packing phase, the main difference between the 2S-1D method and the 1DBP+VRP method is that the 2S-1D uses pallet J4[8,3], while the 1DBP+VRP uses a smaller and cheaper pallet J6[4,2]. (The 1DVRP+BP method is infeasible since the volume of groups of boxes with same destination exceed the maximum pallet capacity.) Consequently, the 2S-1D requires truck K2[20,5], while 1DBP+VRP uses a smaller and cheaper truck K5[12,3]. It leads to a result that the total packing cost of 2S-1D is higher than 1DBP+VRP. However, since the 2S-1D uses a truck with higher capacity, it has an advantage in the routing phase. We could see that the route used by 2S-1D is much more efficient than the route used by 1DBP+VRP.

\subsection{Multi-dimensional}
\label{sub:3D}

Now, we consider the case with multi-dimensional packing. The instances we used are randomly generated as:
\begin{itemize}
    \item ins-7: Set: $|I| = 5$, $|J| = 5$, $|K| = 1$, $|D| = 5$; Parameter: $v \in [1,5]$, $l^0, w^0, h^0 \in [1,5]$, $V \in [5,15]$, $l^1, w^1, h^1 \in [1,6]$, $T \in [36,48]$, $l^2, w^2 \in [5,8]$, $c \in [5,10]$, $C \in [25,30]$ , $c' \in [0,15]$.
    \item ins-8: Set: $|I| = 9$, $|J| = 5$, $|K| = 3$, $|D| = 5$; Parameter: $v \in [1,6]$, $l^0, w^0, h^0 \in [1,5]$, $V \in [10,15]$, $l^1, w^1, h^1 \in [1,6]$, $T \in [20,25]$, $l^2, w^2 \in [2,10]$, $c \in [5,10]$, $C = [25,35]$ , $c' \in [0,15]$.
    \item ins-9: Set: $|I| = 20$, $|J| = 10$, $|K| = 6$, $|D| = 5$; Parameter: $v \in [1,6]$, $l^0, w^0, h^0 \in [1,5]$, $V \in [3,15]$, $l^1, w^1, h^1 \in [1,5]$, $T \in [12,25]$, $l^2, w^2 \in [2,6]$, $c \in [2,5]$, $C = [3,6]$ , $c' \in [0,15]$.
\end{itemize}
Again, we use disjoint method as benchmark. Here, we call the proposed integrated method as 2S-3D, and we use a competing disjoint method called 3DBP+VRP, in which we consider the 3D bin packing in the first phase and consider the vehicle routing problem in the second phase. We also use 3DBP to denote 3D bin packing and use 3DVRP to denote vehicle routing, separately. We do not introduce the 3DVRP+BP method in this case, as it is not computational efficient to consider the grouping options for boxes according their destinations in 3-dimensional packing. The results can be seen in Table \ref{tab:2S-3Dvs3DBP+3DVRP}.

\begin{table}[htb]
\caption{2S-3D vs. 3DBP+VRP cost \& time (in seconds)}
\label{tab:2S-3Dvs3DBP+3DVRP}
\centering
\resizebox{\linewidth}{!}{
\begin{tabular}{ccccccc}
\toprule
\multicolumn{1}{c}{} & \multicolumn{2}{c}{\textbf{Integrated method}} & \multicolumn{4}{c}{\textbf{Disjoint method}}\\
\cmidrule(l{3pt}r{3pt}){2-3} \cmidrule(l{3pt}r{3pt}){4-7}
& 2S-3D cost & 2S-3D time & 3DBP cost & 3DVRP cost  & 3DBP+VRP cost & 3DBP+VRP time\\
\midrule
ins-7 & 65 & 0.15 & 55 & 10 & 65 & 0.36\\
ins-8 & 122 & 0.69 & 102 & 40 & 142 & 1.23\\
ins-9 & 74 & 32.76 & 38 & 70 & 108 & 28.45\\
\bottomrule
\end{tabular}}
\end{table}

\begin{figure}[htb]
\centering
\caption{3D packing structure by 2S-3D of pallet J3 for `ins-9'}
\includegraphics[width=0.65\linewidth]{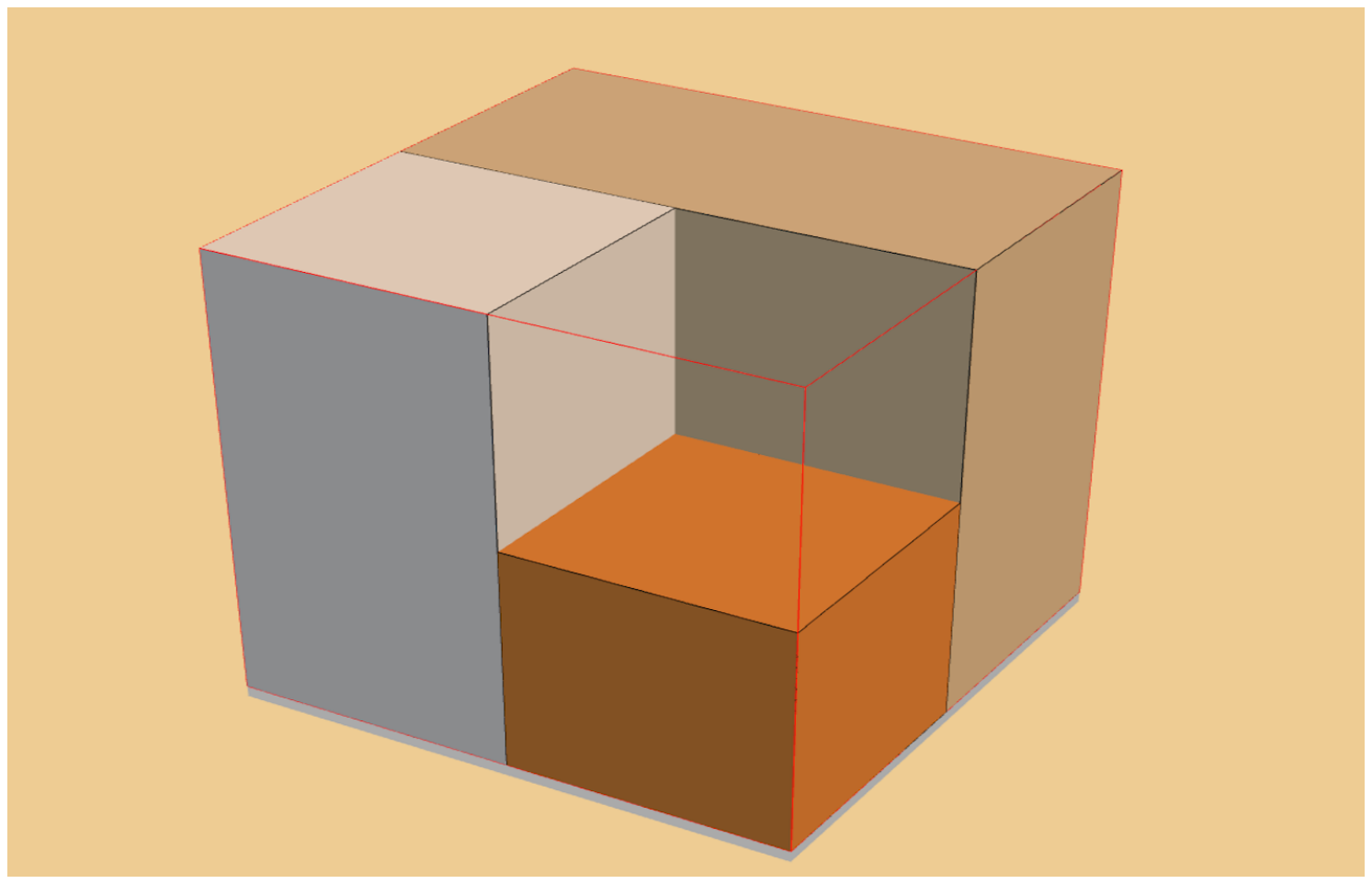}
\label{fig:3D}
\end{figure}

It can be seen from Table \ref{tab:2S-3Dvs3DBP+3DVRP} that the proposed integrated method can again outperform the disjoint method on multi-dimensional packing. Moreover, both methods require similar running time on large-size instances, and the proposed integrated method even runs faster than the disjoint method on small-size instances. We could again breakdown an instance to observe the detailed differences. Here, we use `ins-9' as an example. A 3D packing example at pallet level can be seen in Figure \ref{fig:3D}. Full visualisation results for both methods can be seen in Figure \ref{fig:detail_packing_1} - \ref{fig:detail_packing_2_cont.} in \ref{app:a}. 

\begin{itemize}
    \item 2S-3D: Pallet used: J1[15,5], J2[9,4], J3[8,3], J4[8,3], J5[8,3], J6[4,2], J7[4,2], J9[4,2], J10[3,1]; Truck used: K1[25,6], K2[20,5], K3[12,3], K5[12,3]; Route used: (D0, D2, D4, D0)[11], (D0, D3, D1, D0)[7]. (D0, D1, D0)[4], (D0, D2, D5, D0)[10].
    \item 3DBP+VRP: Pallet used: J1[15,5], J2[9,4], J3[8,3], J4[8,3], J5[8,3], J6[4,2], J8[4,2], J10[3,1]; Truck used: K1[25,6], K3[12,3], K4[12,3], K5[12,3]; Route used: (D0, D5, D4, D2, D1, D3, D0)[19], (D0, D1, D2, D5, D0)[13]. (D0, D4, D1, D0)[21], (D0, D3, D1, D2, D4, D0)[17].
\end{itemize}

It appears that the major difference between 2S-3D and 3DBP+VRP is on the usage of truck K2[20,5] and K4[12,3]. As an integrated method, the 2S-3D prefers to use a larger (but more expensive) truck, since it can cover a wider range of destinations and therefore reduce cost on routing phase. On the other hand, the disjoint 3DBP+VRP method uses a smaller (but cheaper) truck due to its nature of cost minimisation of packing phase. 

\begin{figure}[!ht]
    \centering
    \caption{Packing structure at truck level for `ins-9'}
    \label{fig:truck_packing}
    \begin{subfigure}[b]{0.48\textwidth}
    \centering
    \includegraphics[width=0.85\textwidth]{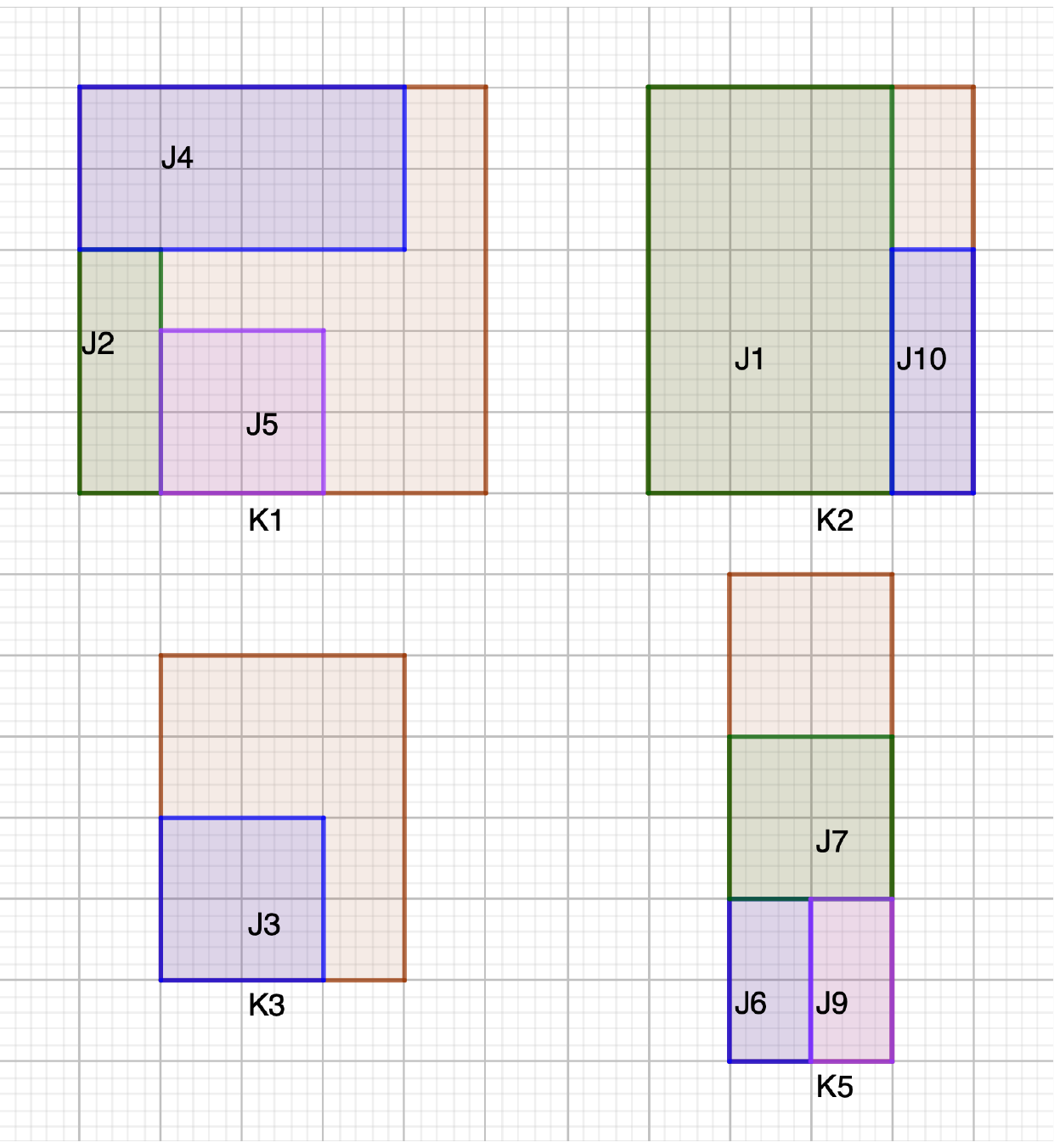}
    \caption{integrated 2S-3D method}
    \end{subfigure}
    \hfill
    \begin{subfigure}[b]{0.48\textwidth}
    \centering
    \includegraphics[width=0.85\textwidth]{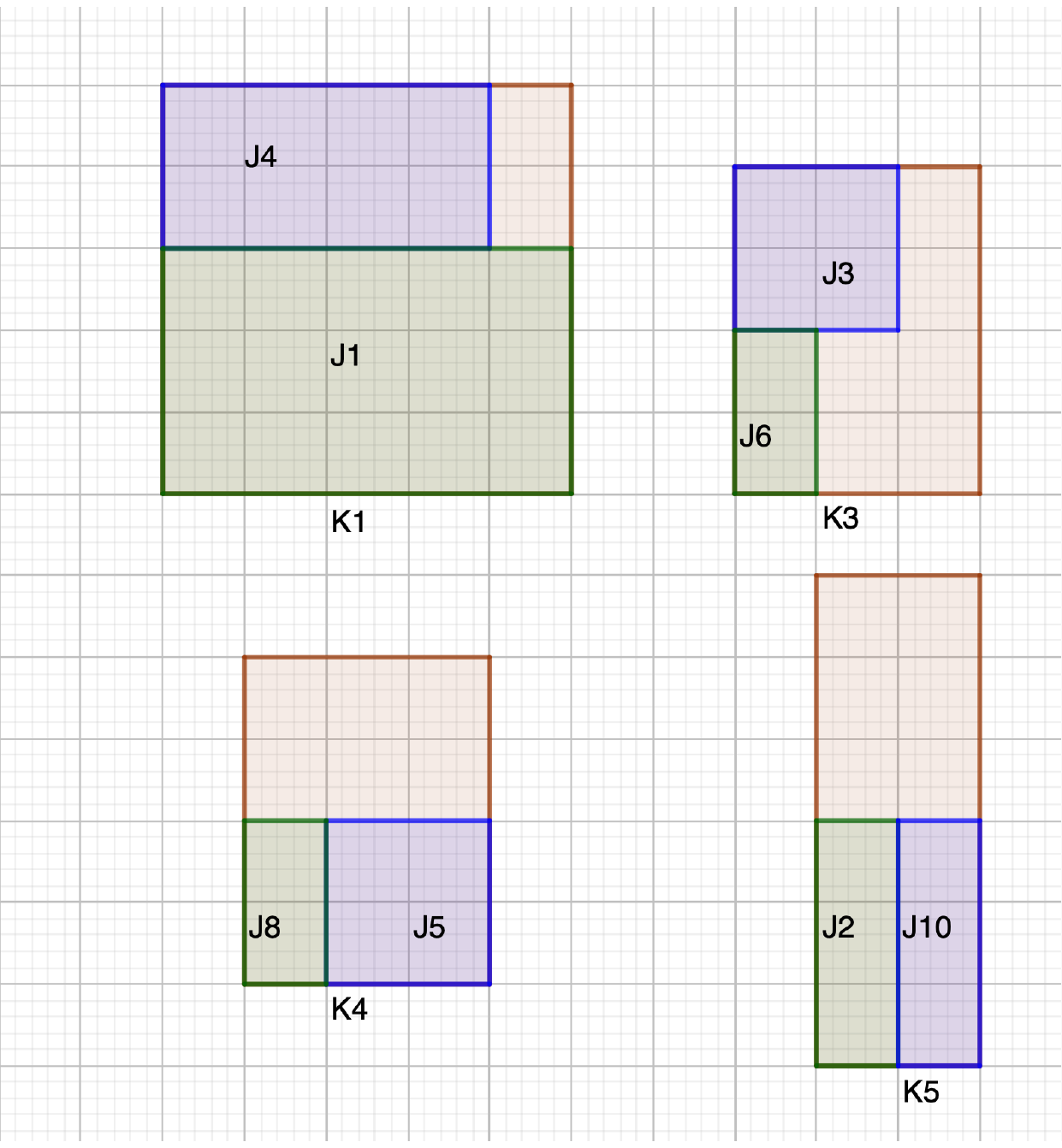}
    \caption{disjoint 3DBP+VRP method}
    \end{subfigure}
\end{figure}

In Figure \ref{fig:truck_packing}, we visualise the packing structures generated by 2S-3D and 3DBP at truck level for `ins-9'. One can see that the packing structure used by these two methods are very different. This is because the integrated 2S-3D method considers the destinations of pallets when packing at truck level, while the disjoint 3DBP+VRP method does not. For instance, the 2S-3D method packs J2, J4 and J5 together in truck K1. In this case, K1 only have to go to 2 destinations, since J2 and J5 need to go to D4, and J4 needs to go to D2. On the other hand, the 3DBP+VRP method packs J1 and J4 in truck K1. It seems like in this case fewer number of pallets are packed in truck K1. However, it needs to go to more destinations (D1, D2, D3, D4 and D5). Therefore, the cost of routing for K1 under 3DBP+VRP method (19) is much higher than that for K1 under 2S-3D method (11). The situations for other trucks are similar to truck K1. Clearly, by not considering the destinations, the disjoint method can achieve lower packing cost, but create much higher routing expense.

\subsection{Real-life instance}
\label{sub:real}

Here, we examine the proposed method again with real-life instance, where the data structure is more complicated. To maintain consistency, we use same solver and same machine as we used in the previous experiments. The measurements and cost data of the boxes/pallets/trucks came from an anonymous data provider. Yet, the destination data is randomly generated as in previous experiment, since it is not provided due to confidentiality. The measurement data can be seen in Table \ref{tab:real-measurement}, and the variable cost of travelling across destinations can be seen in Table \ref{tab:real-travel}.

\begin{table}[htb]
\caption{The measurements of boxes/pallets/trucks for real-life instance}
\label{tab:real-measurement}
\centering
\resizebox{0.58\linewidth}{!}{
\begin{tabular}{cccccc}
\toprule
\multicolumn{1}{c}{} & \multicolumn{5}{c}{\textbf{Boxes}}\\
\cmidrule(l{3pt}r{3pt}){2-6}
& Length & Width & Height & Volume & Destination\\
\midrule
I1 & 475 & 475 & 475 & 11 & D1\\
I2 & 475 & 475 & 475 & 11 & D1\\
I3 & 475 & 475 & 475 & 11 & D2\\
I4 & 326 & 264 & 127 & 1 & D3\\
I5 & 326 & 264 & 127 & 1 & D3\\
I6 & 360 & 350 & 416 & 5 & D4\\
I7 & 490 & 425 & 100 & 2 & D4\\
I8 & 540 & 335 & 195 & 4 & D2\\
I9 & 560 & 325 & 345 & 6 & D5\\
I10 & 312 & 270 & 170 & 1 & D3\\
I11 & 470 & 365 & 415 & 7 & D5\\
I12 & 550 & 340 & 130 & 2 & D5\\
I13 & 380 & 345 & 110 & 1 & D1\\
I14 & 310 & 270 & 130 & 1 & D3\\
I15 & 395 & 165 & 167 & 1 & D1\\
I16 & 395 & 165 & 167 & 1 & D2\\
I17 & 390 & 315 & 175 & 2 & D5\\
I18 & 310 & 240 & 240 & 2 & D4\\
I19 & 525 & 405 & 125 & 3 & D4\\
\midrule
\multicolumn{1}{c}{} & \multicolumn{5}{c}{\textbf{Pallets}}\\
\cmidrule(l{3pt}r{3pt}){2-6}
& Length & Width & Height & Volume & Cost ('00 £)\\
\midrule
J1 & 800 & 800 & 800 & 51 & 6\\
J2 & 800 & 800 & 800 & 51 & 6\\
J3 & 1000 & 800 & 500 & 40 & 5\\
J4 & 1000 & 700 & 500 & 35 & 4\\
J5 & 900 & 400 & 1000 & 36 & 4\\
J6 & 900 & 300 & 1000 & 27 & 2\\
\midrule
\multicolumn{1}{c}{} & \multicolumn{5}{c}{\textbf{Trucks}}\\
\cmidrule(l{3pt}r{3pt}){2-6}
& Length & Width & Height & Volume & Cost ('00 £)\\
\midrule
K1 & 1600 & 1000 & 1000 & 160 & 10\\
K2 & 1000 & 1200 & 800 & 96 & 6\\
K3 & 500 & 1000 & 800 & 40 & 3\\
\bottomrule
\end{tabular}}
\end{table}

\begin{table}[!ht]
\caption{Variable cost of travelling
across destinations ('00 £)}
\label{tab:real-travel}
\centering
\resizebox{0.4\linewidth}{!}{
\begin{tabular}{ccccccc}
\toprule
& D0 & D1 & D2 & D3 & D4 & D5\\
\midrule
D0 & 0 & 2 & 2 & 3 & 8 & 7\\
D1 & 2 & 0 & 3 & 2 & 11 & 13\\
D2 & 2 & 3 & 0 & 5 & 1 & 1\\
D3 & 3 & 2 & 5 & 0 & 14 & 12\\
D4 & 8 & 11 & 1 & 14 & 0 & 3\\
D5 & 7 & 13 & 1 & 12 & 3 & 0\\
\bottomrule
\end{tabular}}
\end{table}

We apply again the proposed integrated 2S-3D method and the competing disjoint 3DBP+VRP method. The detailed results can be seen in Table \ref{tab:real-results}. Here we list the optimal cost, running time and the usage of pallets/trucks. We also give visualisation results for the packing structure of J1 by 2S-3D and J3 by 3DBP+VRP in Figure.

\begin{table}[htb]
\caption{2S-3D vs. 3DBP+VRP on real-life instance}
\label{tab:real-results}
\centering
\resizebox{0.75\linewidth}{!}{
\begin{tabular}{ccc}
\toprule
& \textbf{Integrated 2S-3D method} & \textbf{Disjoint 3DBP+VRP method}\\
\midrule
Optimal cost & 54 & 69\\
Running time & 11.40 & 14.13\\
Pallet in use & J1;J2;J4 & J2;J3;J4\\
Truck in use & K1;K2 & K1;K2\\
Route in use & [D0, D2, D1, D3, D0]; & [D0, D3, D1, D2, D4, D5, D0];\\
& [D0, D2, D4, D5, D0] & [D0, D5, D4, D2, D1, D3, D0]\\
\bottomrule
\end{tabular}}
\end{table}

\begin{figure}[!ht]
    \centering
    \caption{3D packing structure by 2S-3D and 3DBP+VRP for real-life instance}
    \label{fig:real_packing}
    \begin{subfigure}[b]{0.48\textwidth}
    \centering
    \includegraphics[width=0.85\textwidth]{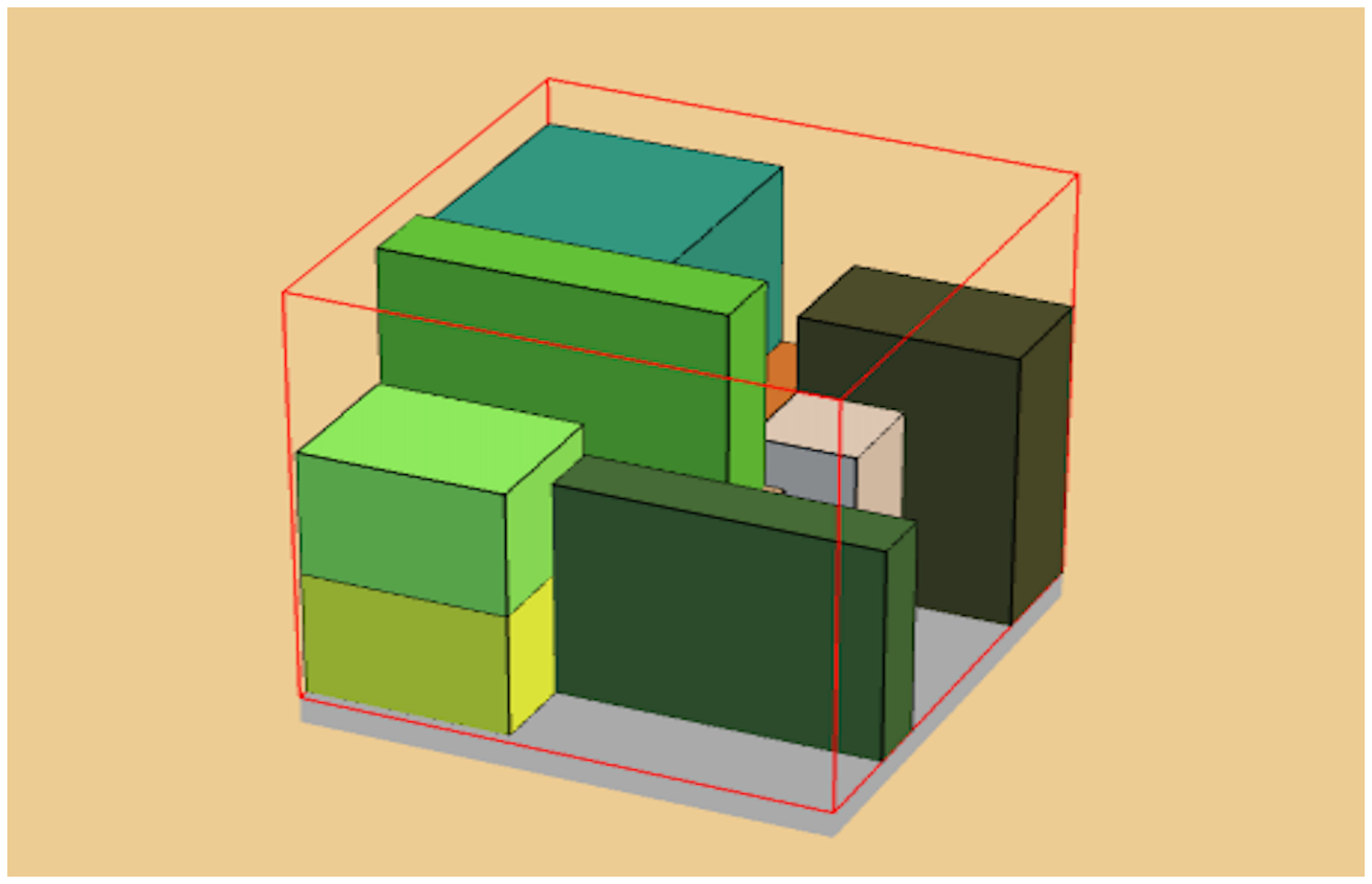}
    \caption{J1 by 2S-3D}
    \end{subfigure}
    \hfill
    \begin{subfigure}[b]{0.48\textwidth}
    \centering
    \includegraphics[width=0.85\textwidth]{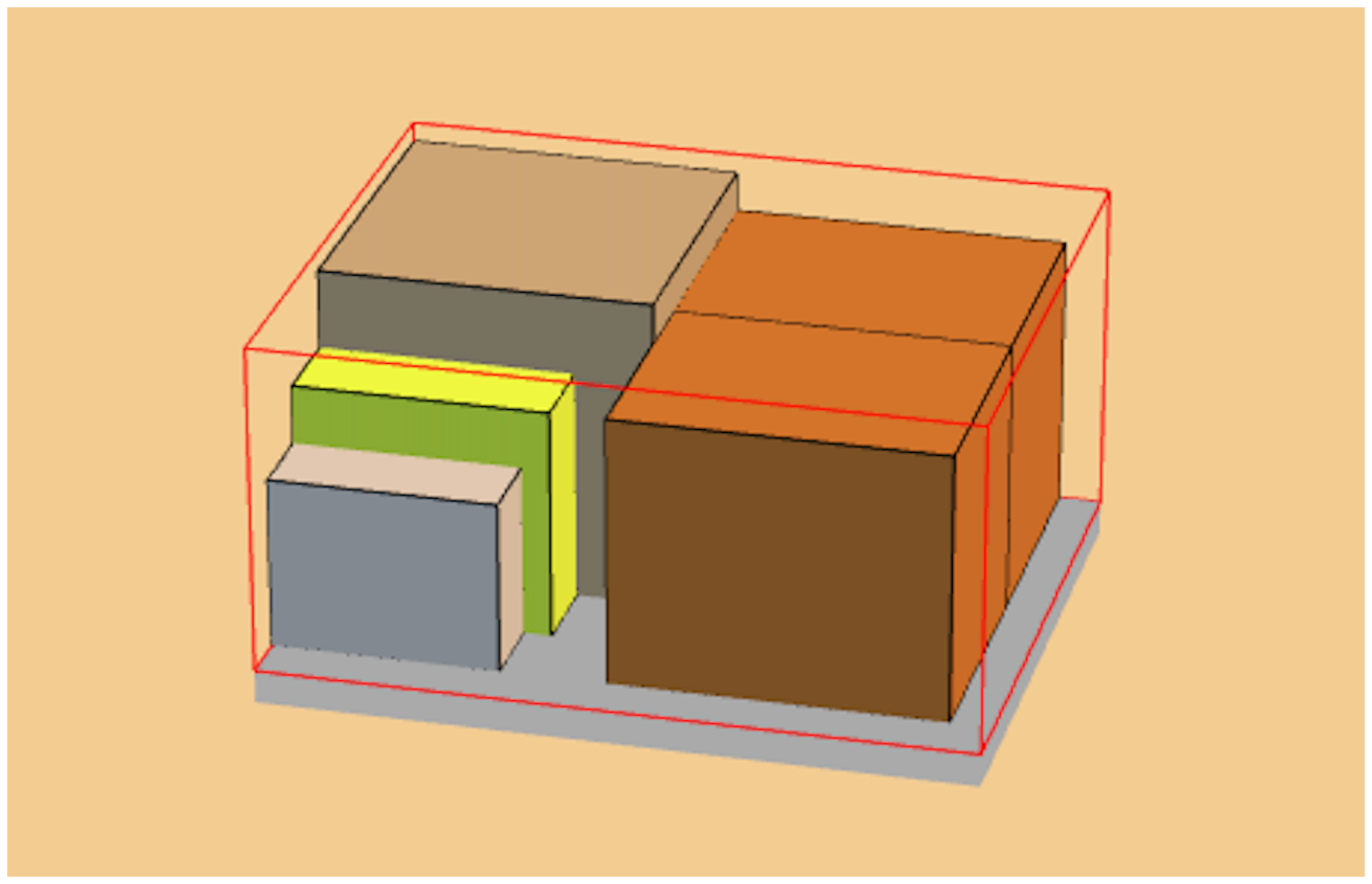}
    \caption{J3 by 3DBP+VRP}
    \end{subfigure}
\end{figure}

The results from real-life instance confirm our findings again. The proposed integrated 2S-3D method and the traditional disjoint method use different packing structures. Although the 2S-3D method uses a more costly pallet J1 rather than J3, it achieves overall lower cost than the disjoint method, as it acquires more efficient routes. Moreover, it does not require extensive computational effect.


\section{Concluding remarks}
\label{se:conclusion}

This paper considered the two-level capacitated vehicle routing problem (2S-CVRP). The idea of the 2S-CVRP is to combine the two-level BPP and the capacitated VRP into an integrated framework. A solution method based on mixed integer program (MIP) was provided in the paper. We then conducted extensive experiment to test the validity of the method. Results show that the proposed integrated approach can successfully capture the nature of the problem and can achieve lower cost than traditional disjoint approach.

An interesting (and challenging) topic for future research is the development of heuristic method for the 2S-CVRP. As an exact method, the MIP could become computational inefficient as the instance size grows. In that case, it could be very useful to have a heuristic method to work with.


\section*{Acknowledgements}
I thank Malek Almousa, Ke Fang and Jing Lyu for illuminating discussions and feedback on the manuscript. I acknowledge financial support from Decision Lab Ltd. 


\bibliographystyle{plain}
\bibliography{2s-cvrp}


\newpage
\appendix

\section{Visualisation results by 2S-3D for `ins-9'}
\label{app:a}

\begin{figure}[!ht]
    \centering
    \caption{3D packing structure by 2S-3D for ‘ins-9’}
    \label{fig:detail_packing_1}
    \begin{subfigure}[b]{0.48\textwidth}
    \centering
    \includegraphics[width=0.8\textwidth]{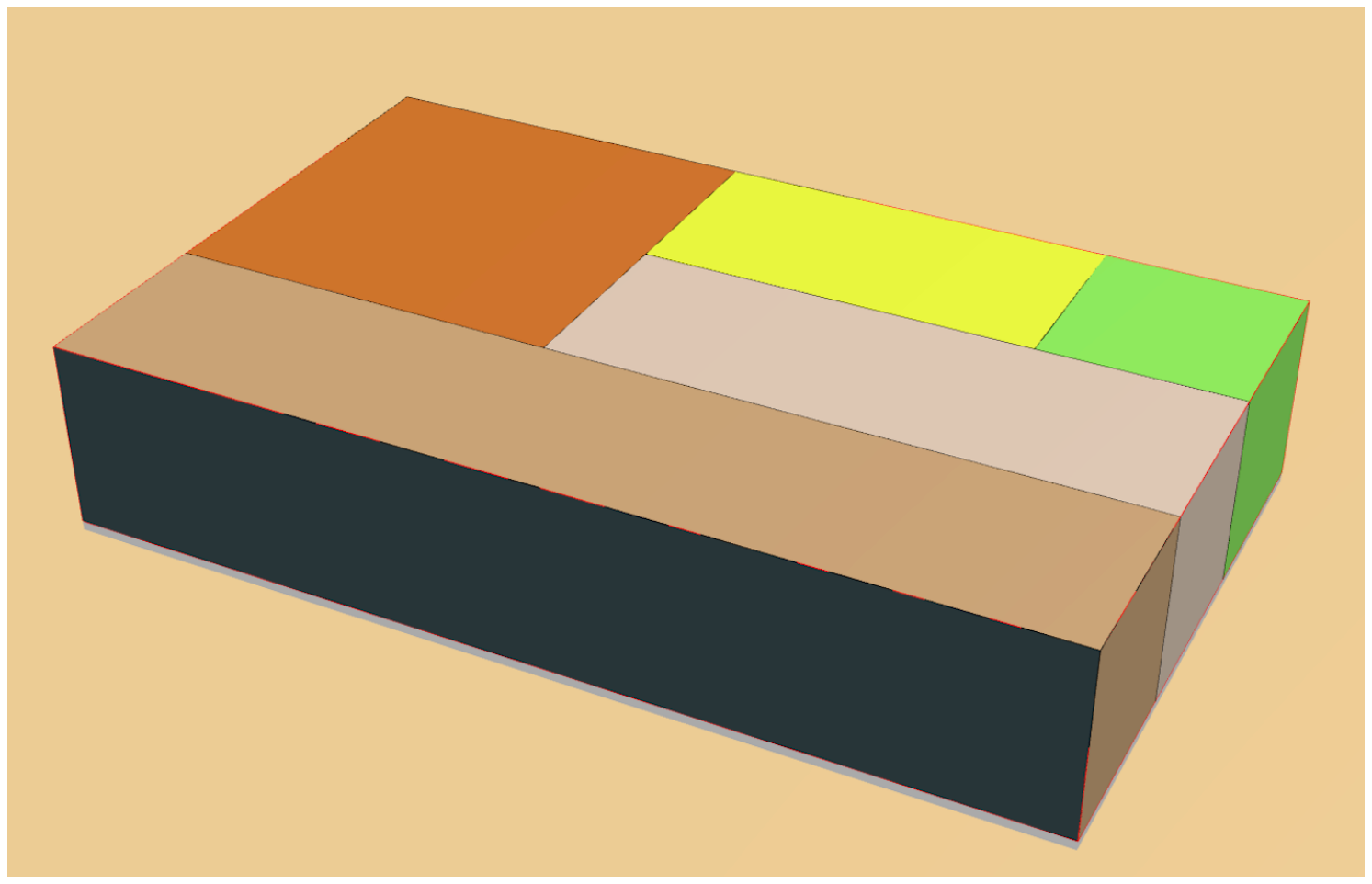}
    \caption{J1}
    \end{subfigure}
    \hfill
    \begin{subfigure}[b]{0.48\textwidth}
    \centering
    \includegraphics[width=0.8\textwidth]{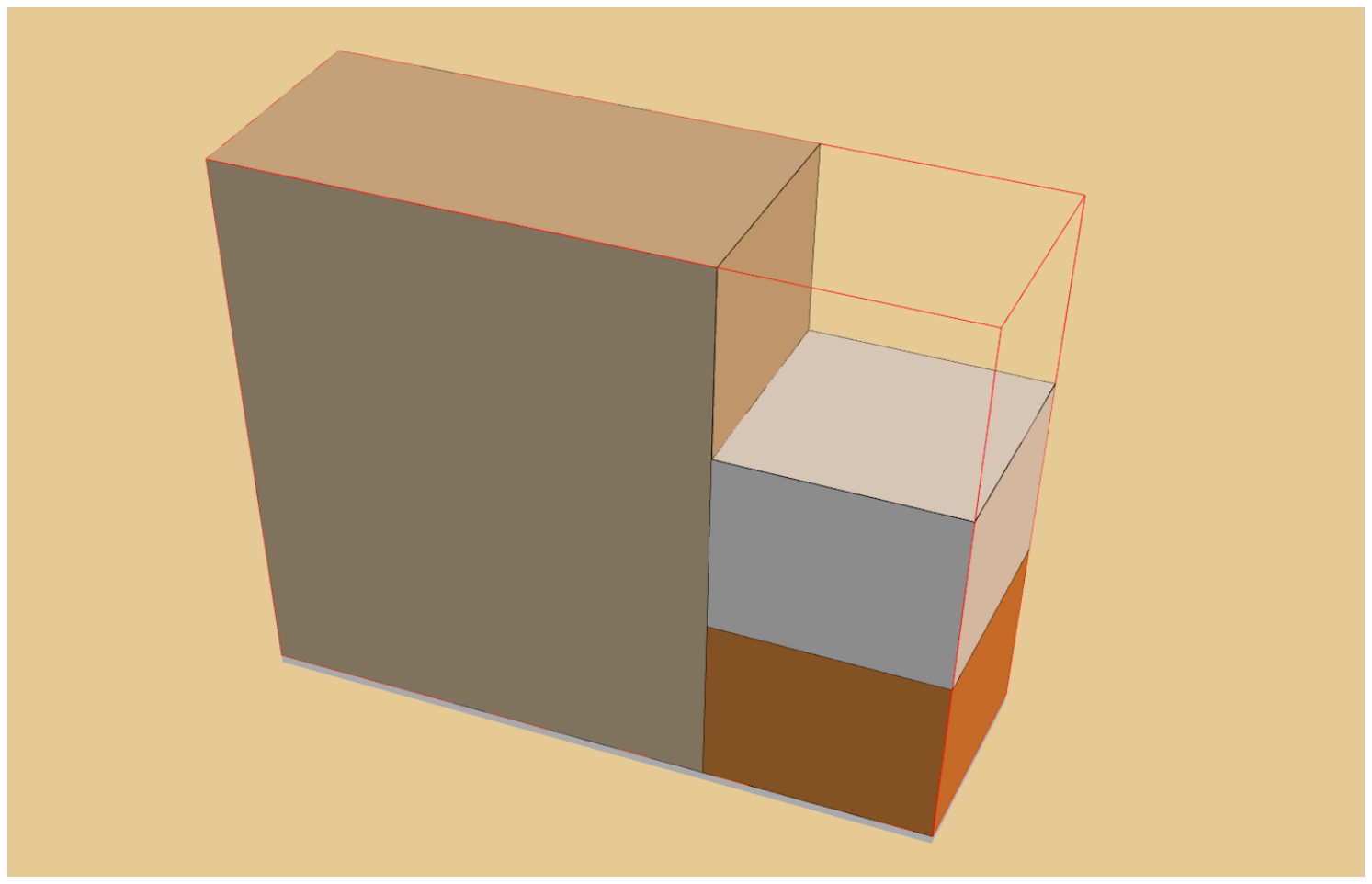}
    \caption{J2}
    \end{subfigure}\\
    \begin{subfigure}[b]{0.48\textwidth}
    \centering
    \includegraphics[width=0.8\textwidth]{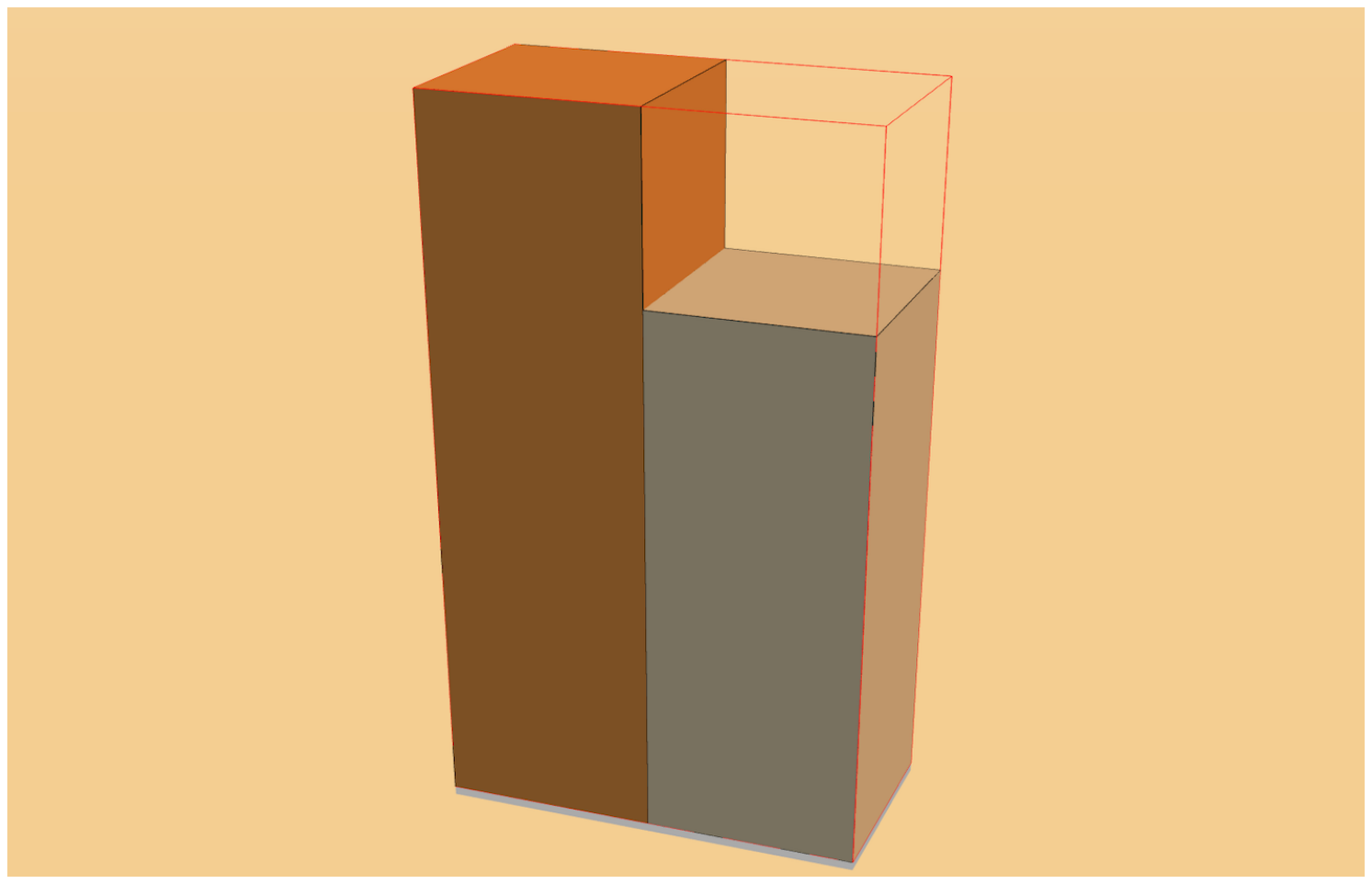}
    \caption{J4}
    \end{subfigure}
    \hfill
    \begin{subfigure}[b]{0.48\textwidth}
    \centering
    \includegraphics[width=0.8\textwidth]{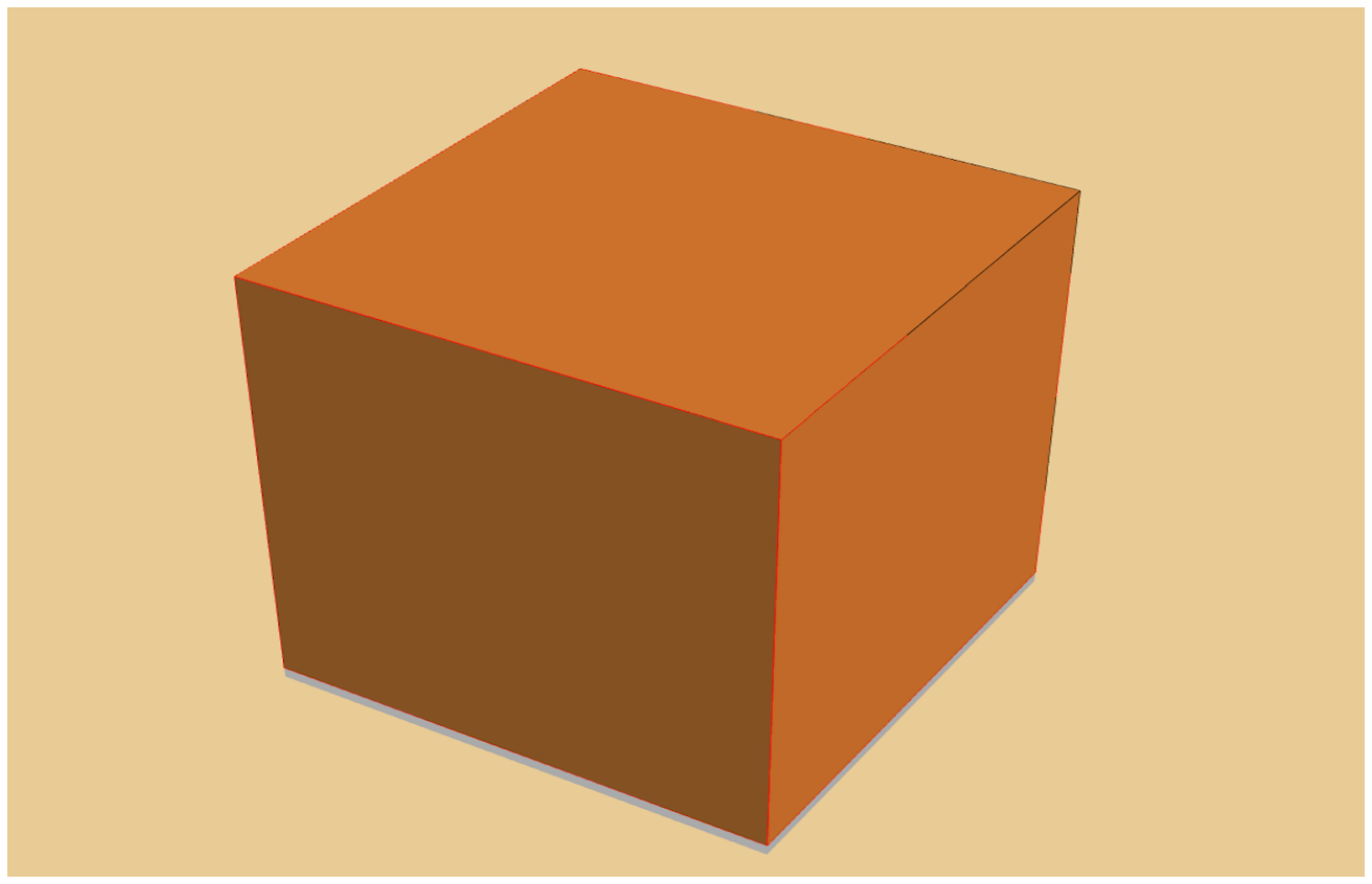}
    \caption{J5}
    \end{subfigure}
\end{figure}

\begin{figure}[!ht]
    \centering
    \caption{3D packing structure by 2S-3D for ‘ins-9’ cont.}
    \label{fig:detail_packing_1_cont.}
    \begin{subfigure}[b]{0.48\textwidth}
    \centering
    \includegraphics[width=0.8\textwidth]{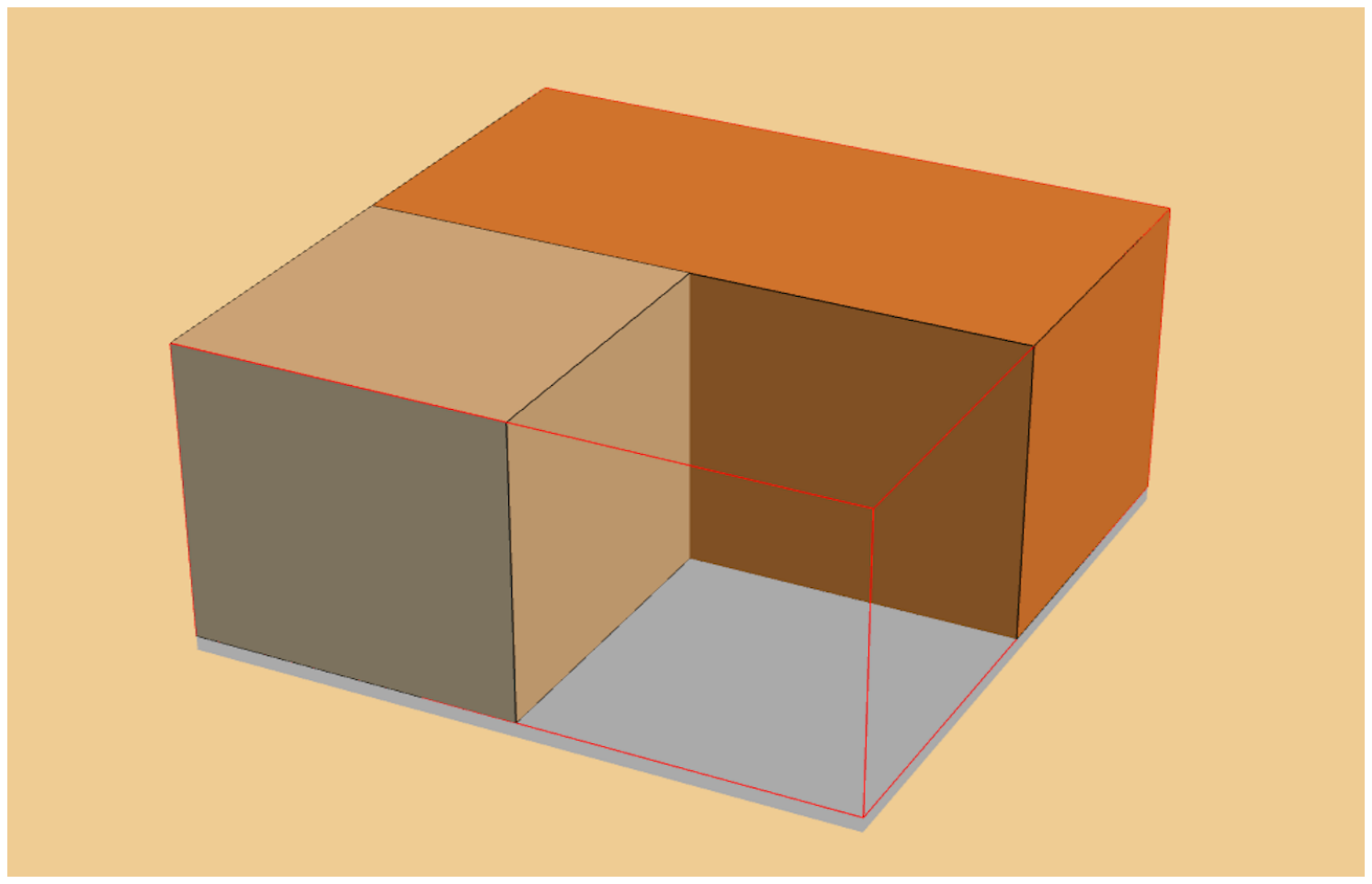}
    \caption{J6}
    \end{subfigure}
    \hfill
    \begin{subfigure}[b]{0.48\textwidth}
    \centering
    \includegraphics[width=0.8\textwidth]{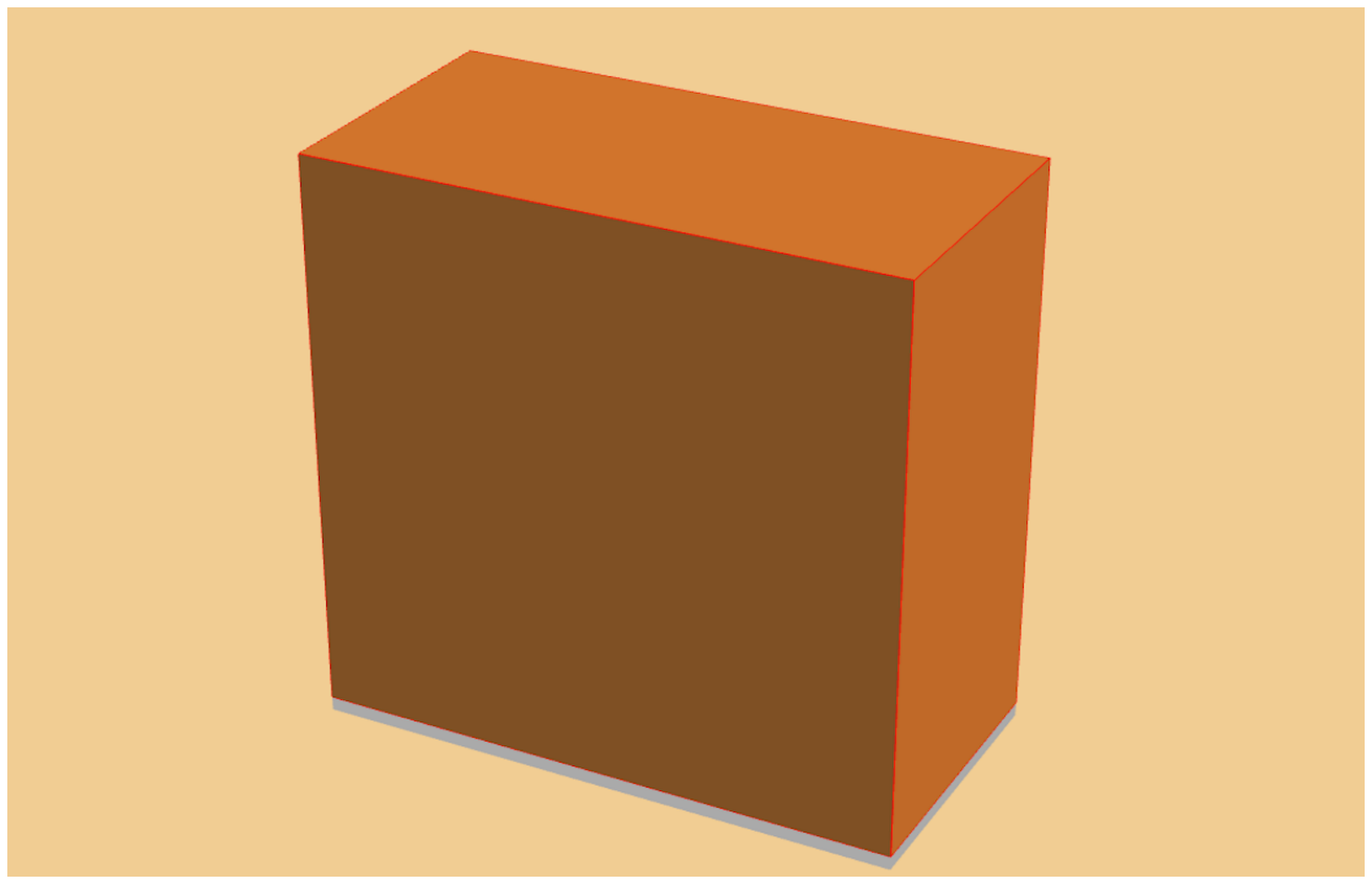}
    \caption{J7}
    \end{subfigure}\\
    \begin{subfigure}[b]{0.48\textwidth}
    \centering
    \includegraphics[width=0.8\textwidth]{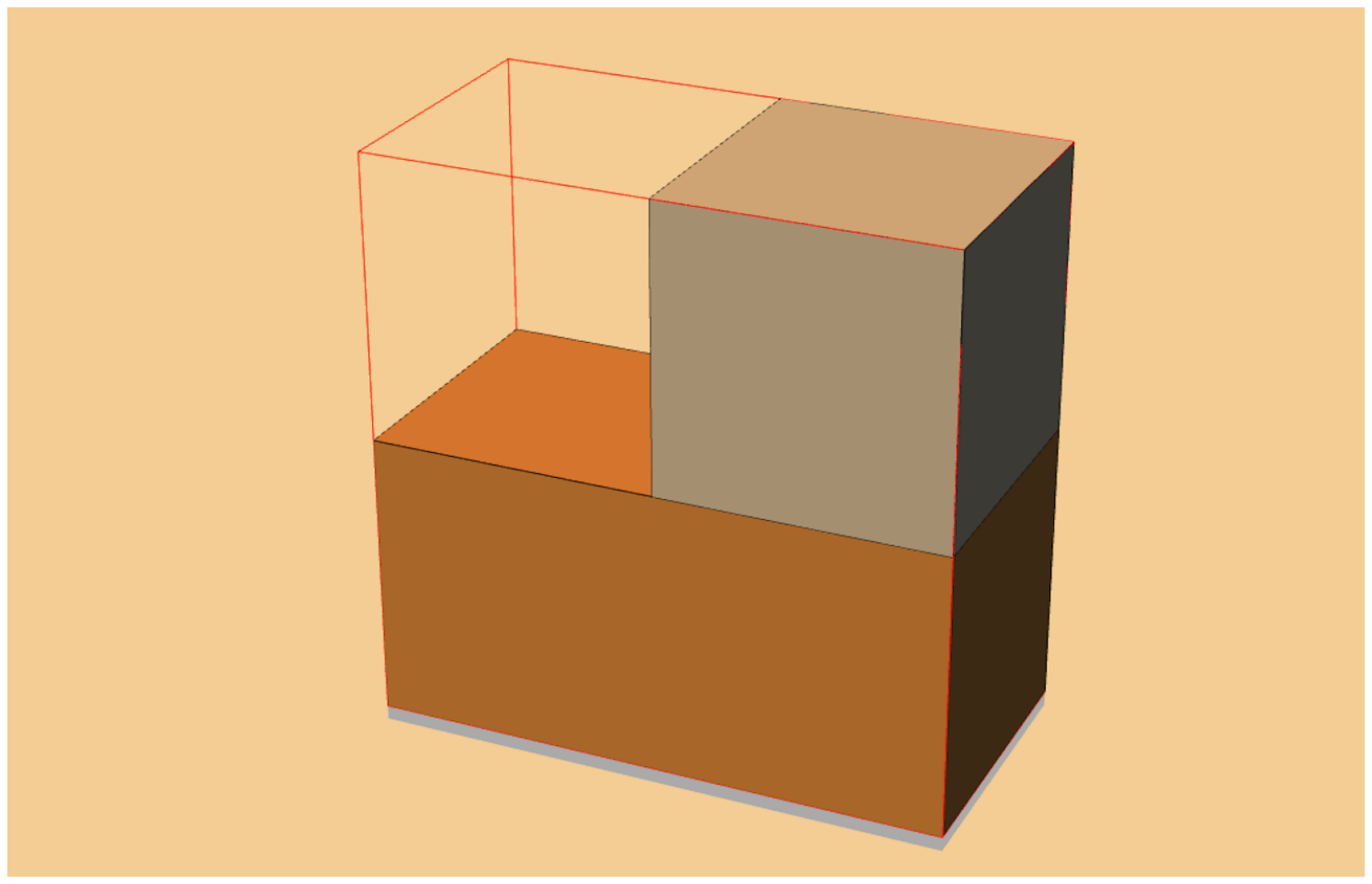}
    \caption{J9}
    \end{subfigure}
    \hfill
    \begin{subfigure}[b]{0.48\textwidth}
    \centering
    \includegraphics[width=0.8\textwidth]{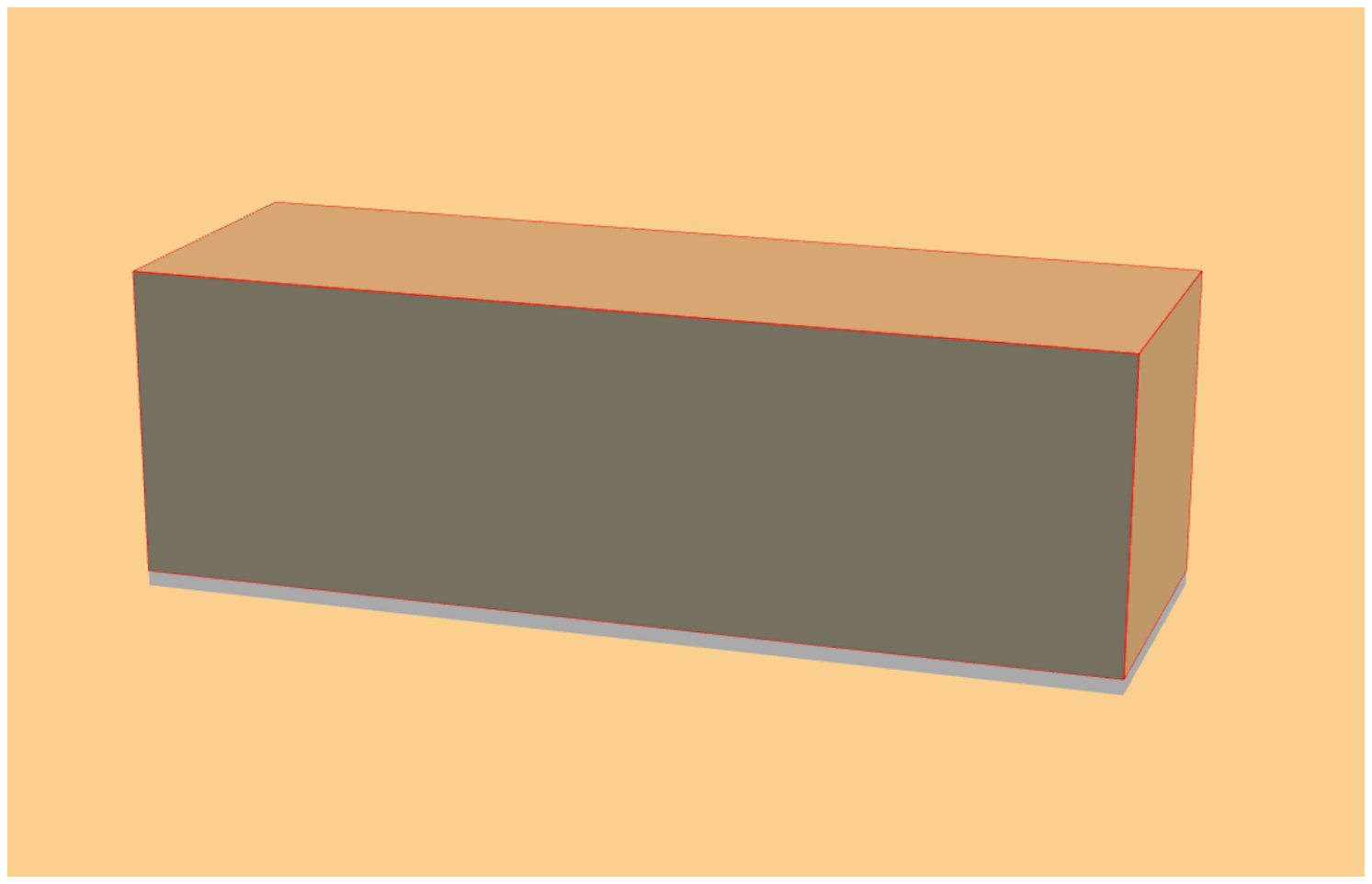}
    \caption{J10}
    \end{subfigure}
\end{figure}

\begin{figure}[!ht]
    \centering
    \caption{3D packing structure by 3DBP for ‘ins-9’}
    \label{fig:detail_packing_2}
    \begin{subfigure}[b]{0.48\textwidth}
    \centering
    \includegraphics[width=0.8\textwidth]{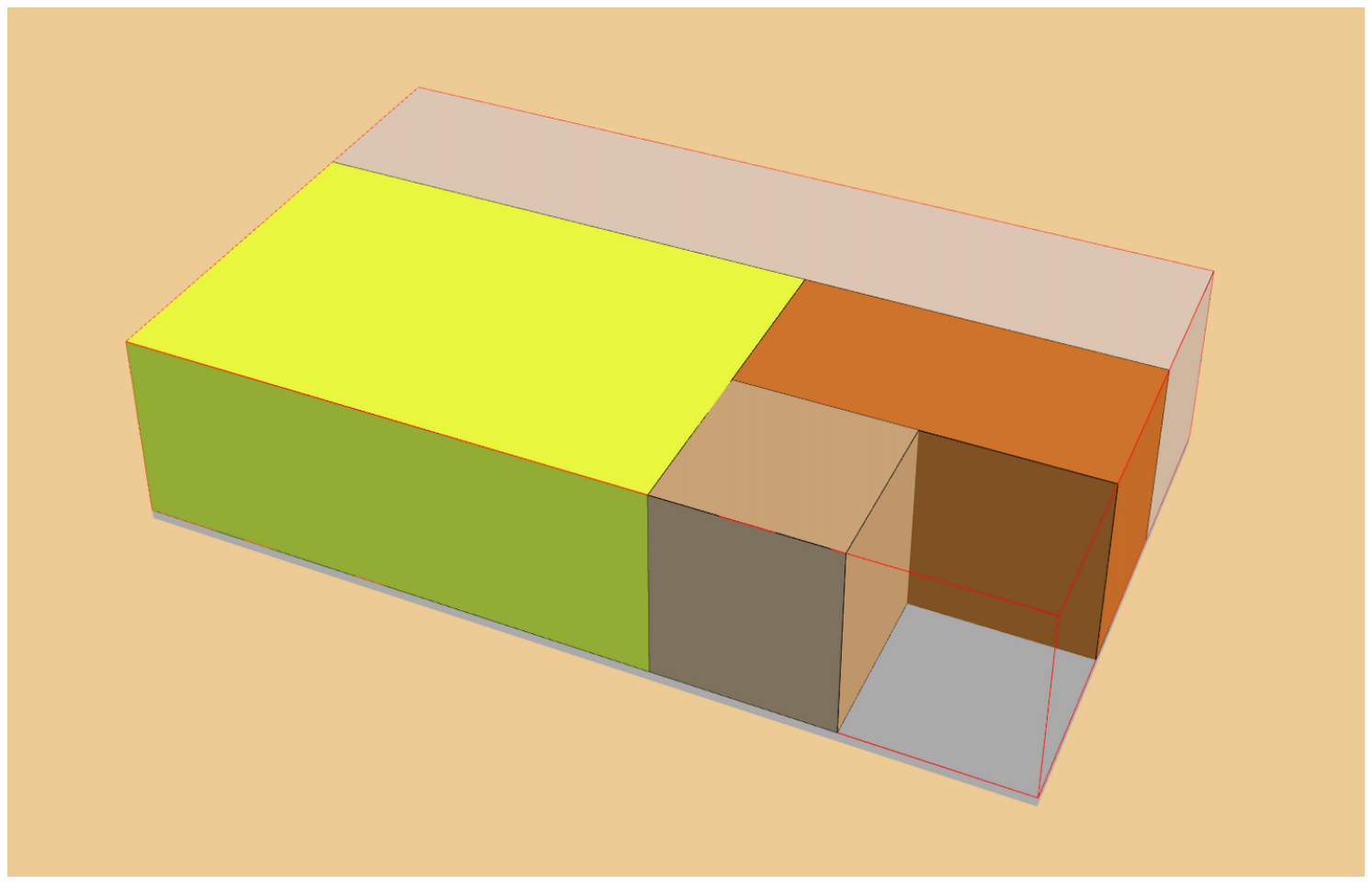}
    \caption{J1}
    \end{subfigure}
    \hfill
    \begin{subfigure}[b]{0.48\textwidth}
    \centering
    \includegraphics[width=0.8\textwidth]{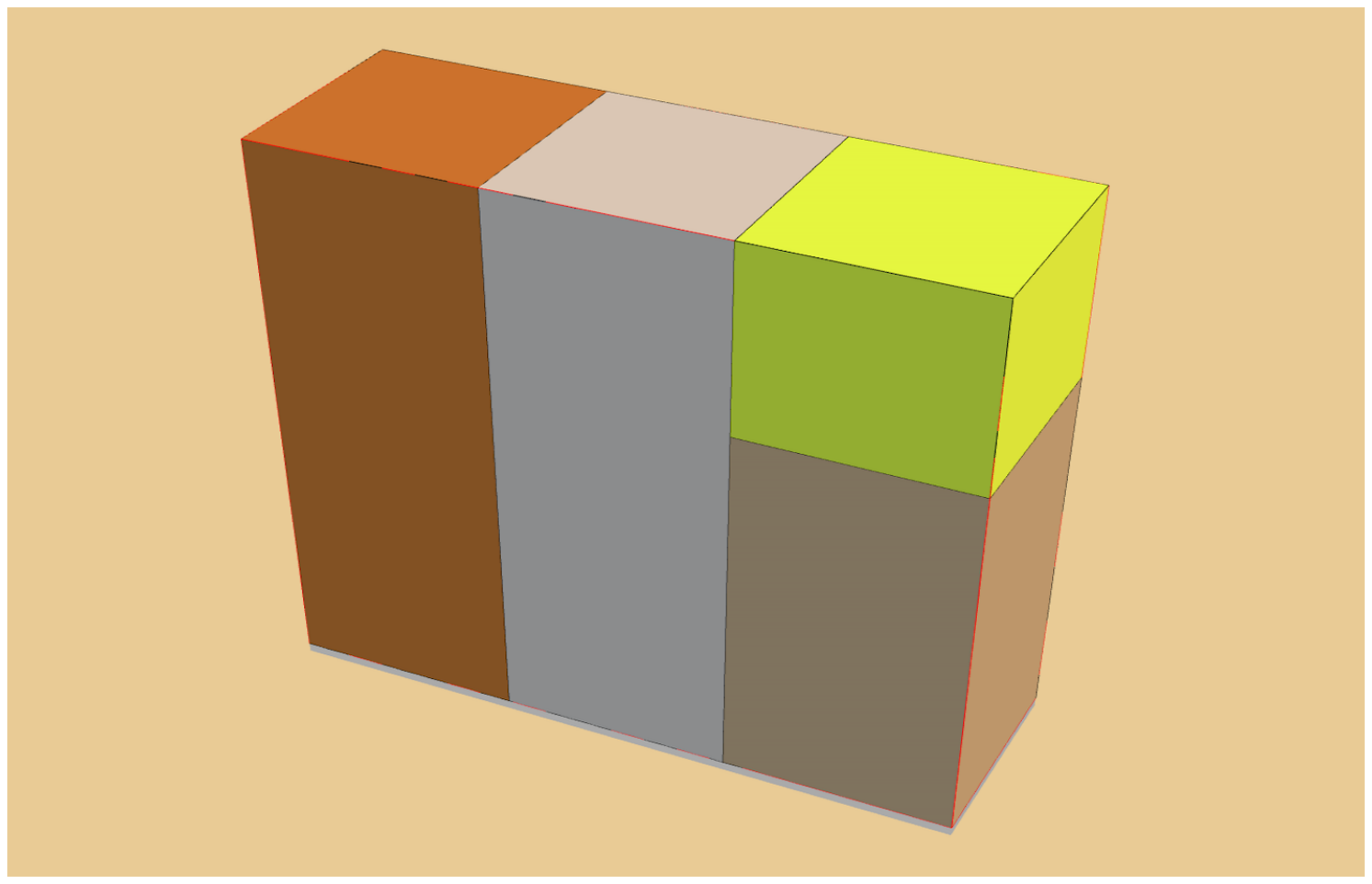}
    \caption{J2}
    \end{subfigure}\\
    \begin{subfigure}[b]{0.48\textwidth}
    \centering
    \includegraphics[width=0.8\textwidth]{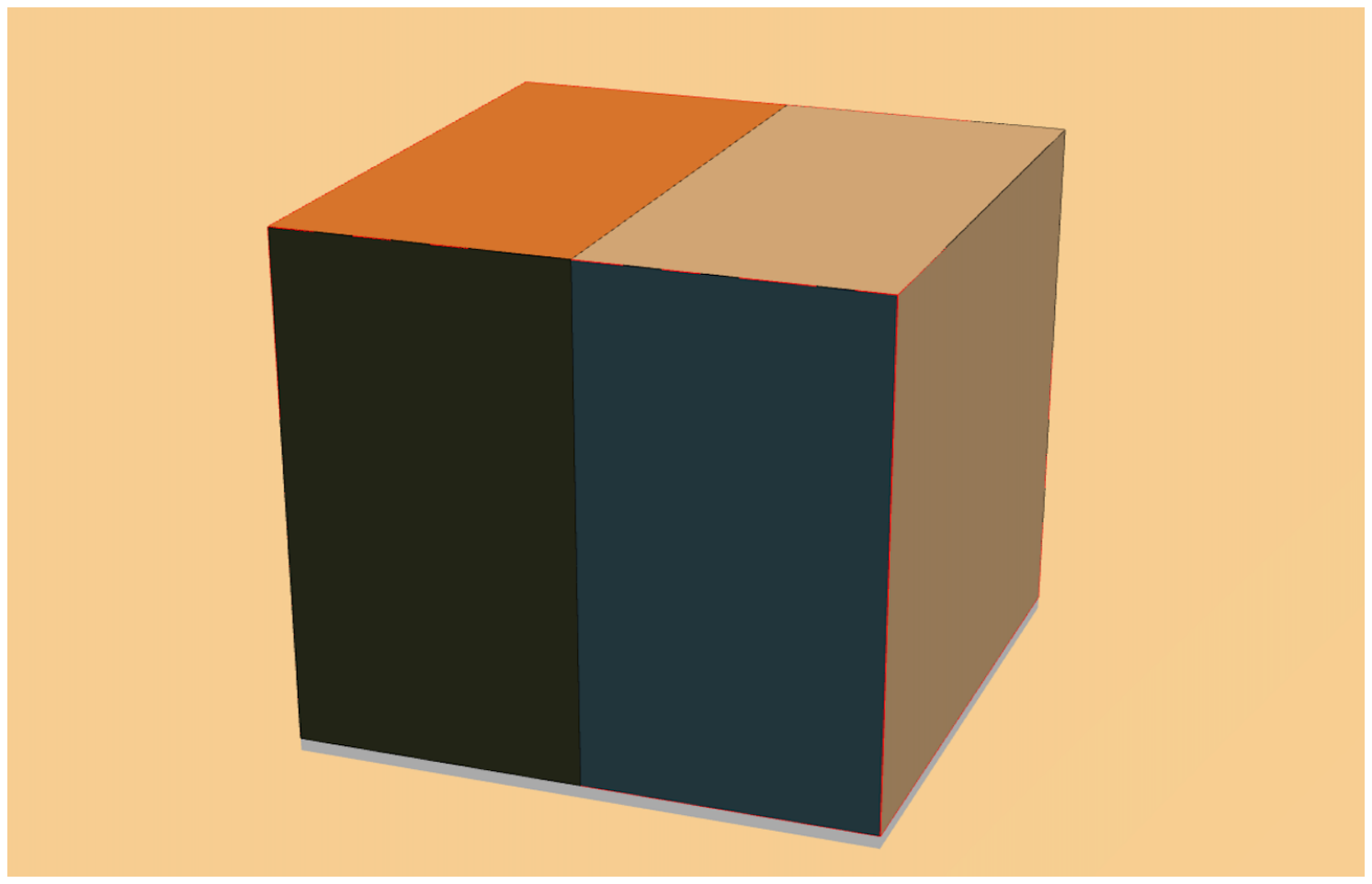}
    \caption{J3}
    \end{subfigure}
    \hfill
    \begin{subfigure}[b]{0.48\textwidth}
    \centering
    \includegraphics[width=0.8\textwidth]{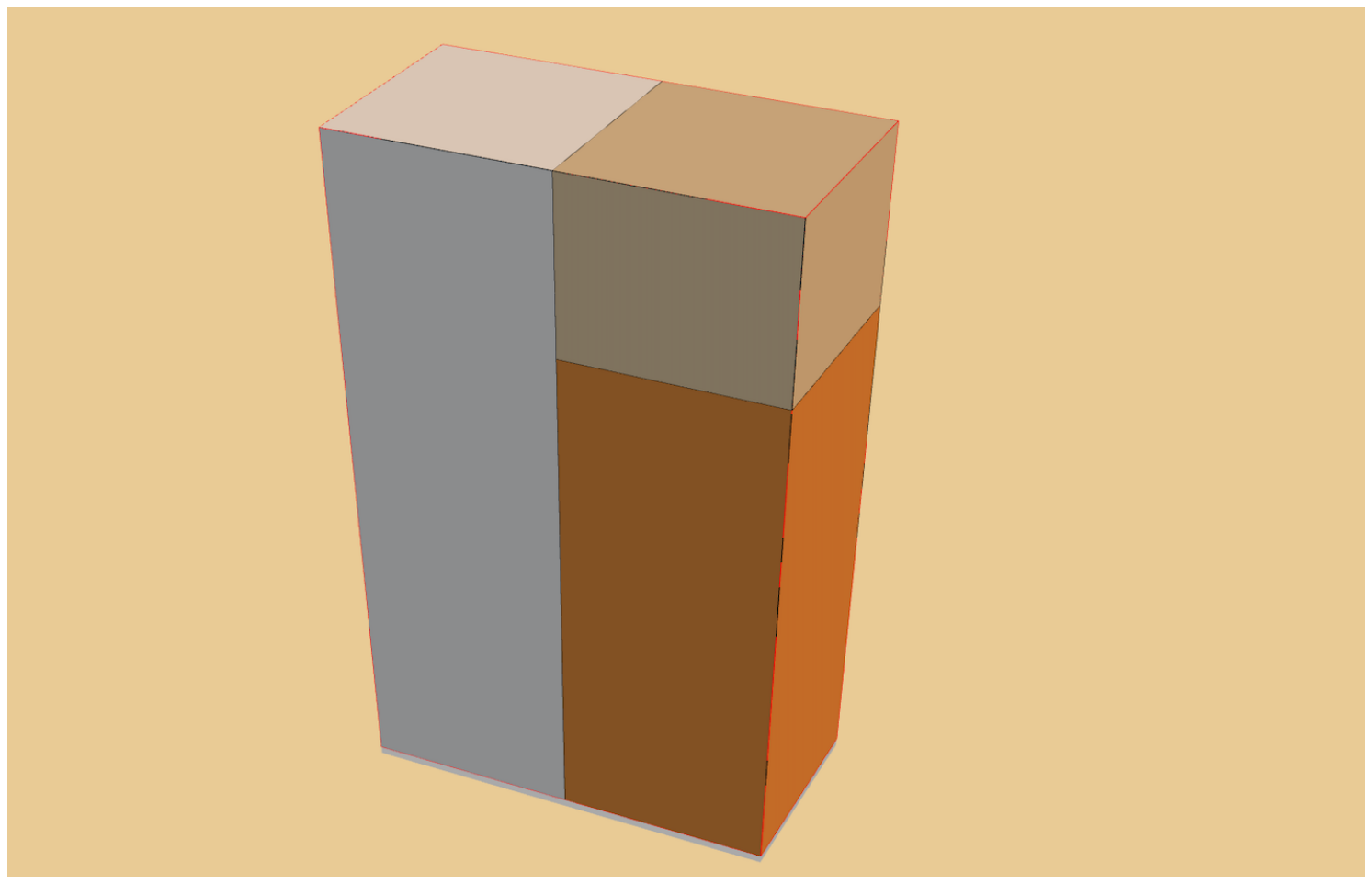}
    \caption{J4}
    \end{subfigure}
\end{figure}

\begin{figure}[!ht]
    \centering
    \caption{3D packing structure by 3DBP for ‘ins-9’ cont.}
    \label{fig:detail_packing_2_cont.}
    \begin{subfigure}[b]{0.48\textwidth}
    \centering
    \includegraphics[width=0.8\textwidth]{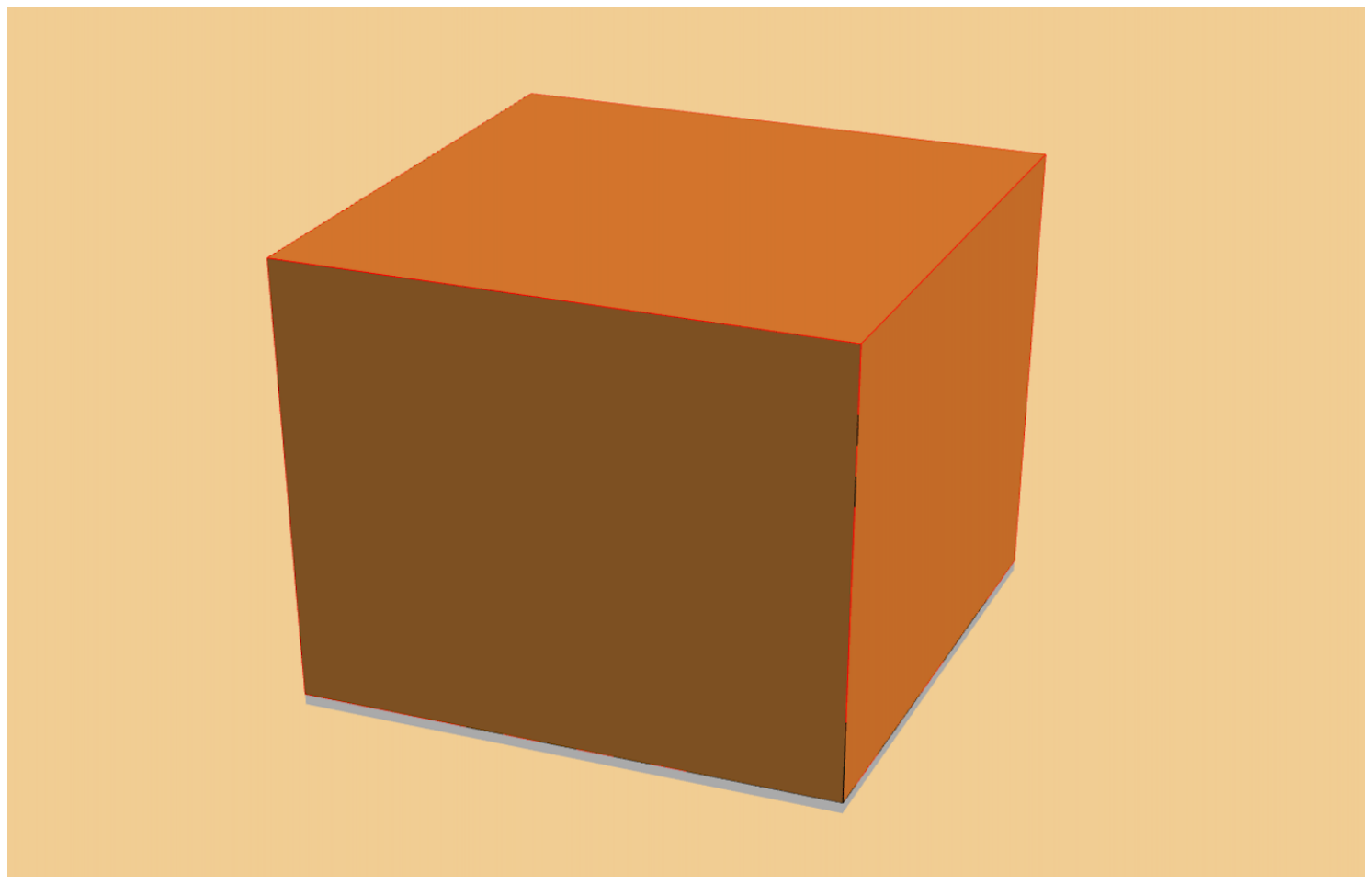}
    \caption{J5}
    \end{subfigure}
    \hfill
    \begin{subfigure}[b]{0.48\textwidth}
    \centering
    \includegraphics[width=0.8\textwidth]{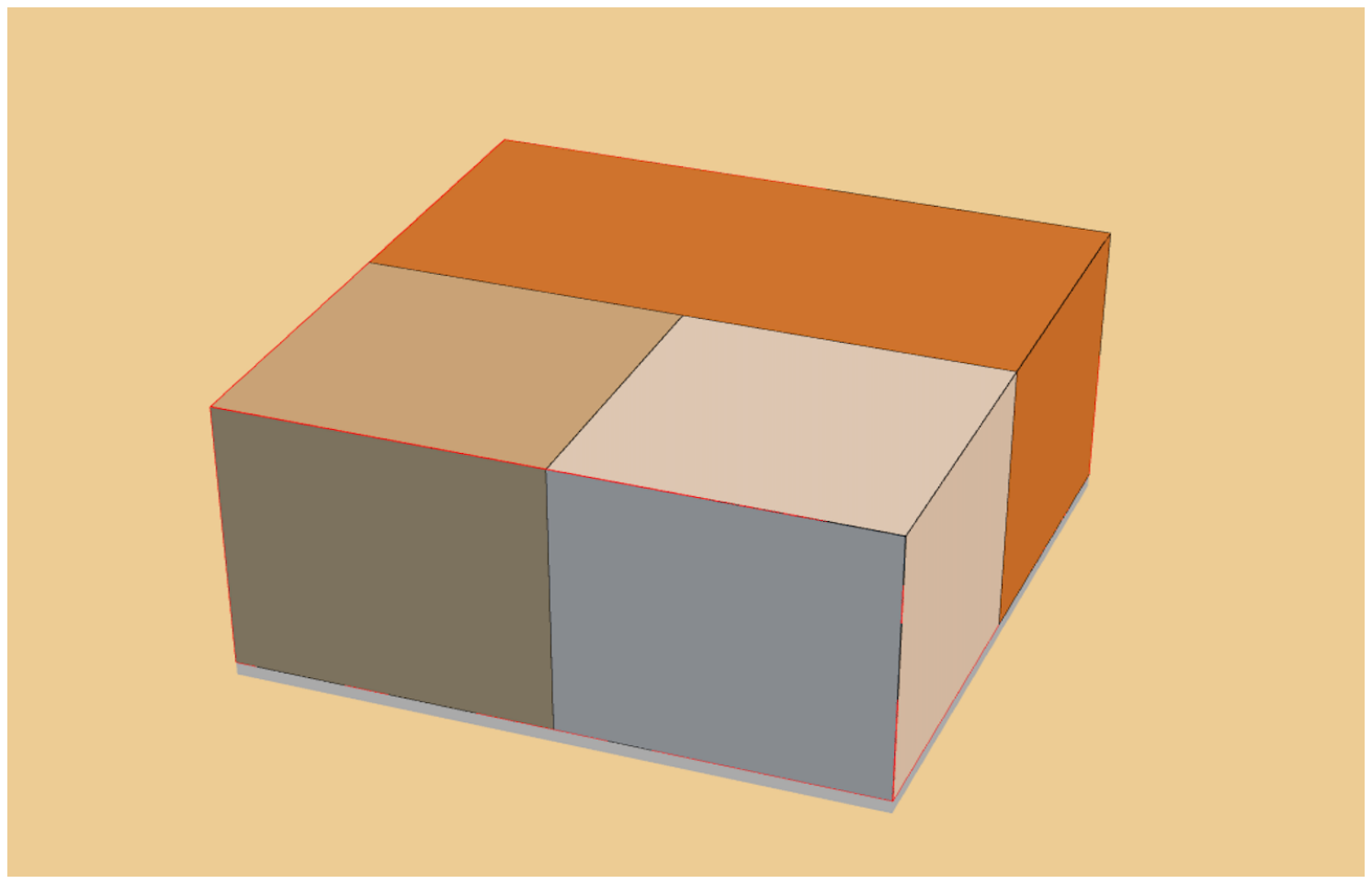}
    \caption{J6}
    \end{subfigure}\\
    \begin{subfigure}[b]{0.48\textwidth}
    \centering
    \includegraphics[width=0.8\textwidth]{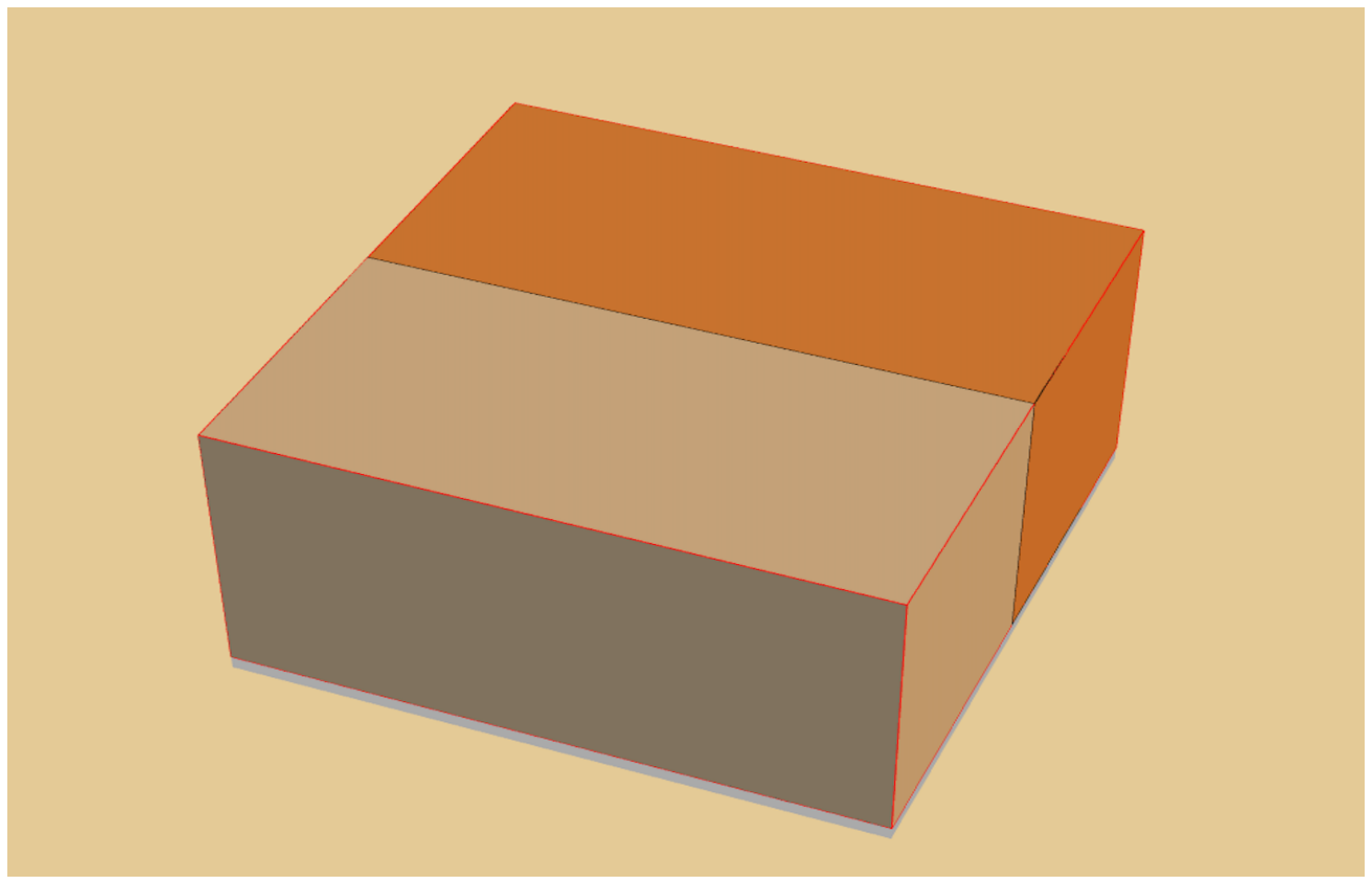}
    \caption{J8}
    \end{subfigure}
    \hfill
    \begin{subfigure}[b]{0.48\textwidth}
    \centering
    \includegraphics[width=0.8\textwidth]{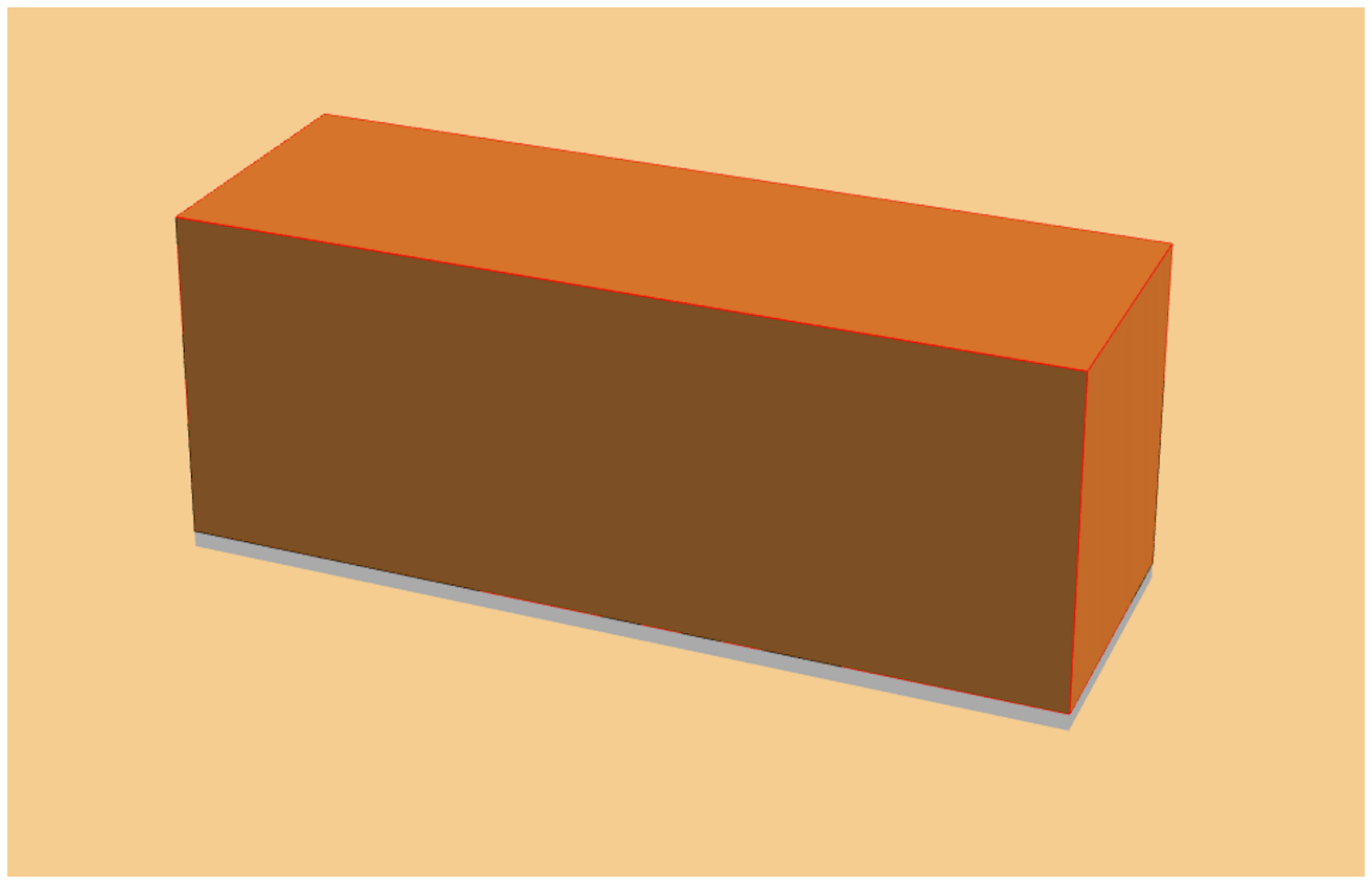}
    \caption{J10}
    \end{subfigure}
\end{figure}

\end{document}